\documentclass[journal]{IEEEtran}
\usepackage{graphicx}
\usepackage{amsfonts}
\usepackage{amsthm}
\usepackage[cmex10]{amsmath}
\usepackage{amssymb,times}
\usepackage{algpseudocode}
\usepackage{algorithm}
\usepackage{url}
\usepackage{array}
\usepackage{subfigure} 
\usepackage{mathtools}
\usepackage{bm,array}
\usepackage{tabularx}
\usepackage{epsfig} 
\usepackage{epstopdf}
\usepackage{color} 
\usepackage[numbers,sort&compress]{natbib}
\usepackage{multirow,booktabs}
\usepackage{alltt}
\usepackage{bigstrut}
\usepackage[autostyle]{csquotes}
\usepackage{times} 
\usepackage{tabularx,ragged2e,booktabs,caption,caption}  
\usepackage{subfig}
\usepackage{balance}
\usepackage{verbatim,enumerate}
\theoremstyle{plain}
\newtheorem{thm}{Theorem} 
\theoremstyle{remark}
\newtheorem{rem}{\bf Remark} 
 
\def \St {{\mathcal{S}_T}} 
\def\Gt {{\mathcal{G}_T}} 
\def \x {{\mathbf{x}}}
\def \y {{\mathbf{y}}}
\def \xh {{\hat{\x}}}
\def \d {{\mathbf{d}}}
\def \v {{\mathbf{v}}}
 
\def \R {{\mathbb{R}}}
\def \N {{\mathbb{N}}}
\def \X {{\mathcal{X}}} 
\def\gt {{\gamma(t)}}
\def\s {{\sum_{t=1}^T}}
\def \A {\mathbf{A}}

\def \bell	{{\boldsymbol{\ell}}}
\def\w{\mathbf{w}}
\def \g {{\mathbf{g}}}
\def \px {{\mathcal{P}_{\X}}}
\def \gmax {{\bar{\gamma}}}
\def \gmin {{\underline{\gamma}}}
\def \alfmax {{\bar{\alpha}}}
\def \alfmin {{\underline{\alpha}}}
\def \Utt {{\tilde{U}_t}}
\def \nt {{\mathbf{n}(t)}} 
\def \E {{\mathop{{}\mathbb{E}}}}
\def \gradt {{\tilde{\nabla}}}
\def \nt {{\mathbf{n}(t)}}

\DeclarePairedDelimiter\abs{\lvert}{\rvert}
\DeclareMathAlphabet{\mathcal}{OMS}{cmsy}{m}{n} 
\newcommand{\norm}[1]{\left\lVert#1\right\rVert}
\renewcommand{\O}[1]{\mathcal{O}\left(#1\right)}
\providecommand{\ip}[1]{\langle#1\rangle}
\makeatletter
\newcommand{\colr}[1]{\textcolor{red}{#1}}
\newcommand{\colb}[1]{\textcolor{black}{#1}}
\makeatother

\begin{document}
\title{{\colb{Online Trajectory Optimization Using Inexact Gradient Feedback for Time-Varying Environments}}}
\author{Mohan Krishna Nutalapati\IEEEauthorrefmark{1}, Amrit Singh Bedi\IEEEauthorrefmark{2}, Ketan Rajawat\IEEEauthorrefmark{1}, and Marceau Coupechoux\IEEEauthorrefmark{3}
\thanks{\IEEEauthorblockA{\IEEEauthorrefmark{1}Department of Electrical
			Engineering (EE), Indian Institute of Technology Kanpur, Kanpur, 208016, India,
			(e-mail: $\left\{\text{nmohank, ketan} \right\}$@iitk.ac.in).} \IEEEauthorblockA{\IEEEauthorrefmark{2}was a Ph.D. student with the Dept. of EE, IIT Kanpur and is currently a postdoc with the US Army Research Laboratory, Adelphi, Maryland, USA, (e-mail: amrit0714@gmail.com).} \IEEEauthorblockA{\IEEEauthorrefmark{3}LTCI, Telecom Paris, IP
			Paris, 91120 Palaiseau France, Paris, France, (email: marceau.coupechoux@telecom-paris.fr).}} }

\IEEEtitleabstractindextext{
\begin{abstract} \colb{This paper considers the problem of online trajectory design under time-varying environments. We formulate the general trajectory optimization problem within the framework of time-varying constrained convex optimization and propose a novel version of online gradient ascent algorithm for such problems. Moreover, the gradient feedback is noisy and allows to use the proposed algorithm for a range of practical applications where it is difficult to acquire the true gradient. Since we are interested in constrained online convex optimization, we carefully select the step size at each iteration so that the iterates stay feasible. In contrast to the most available literature, we present the offline sublinear regret of the proposed algorithm up to the path length variations of the optimal offline solution, the cumulative gradient, and the error in the gradient variations. Furthermore, we establish a lower bound on the offline dynamic regret, which defines the optimality of any trajectory. } 
		
\colb{To show the efficacy of the proposed algorithm, we consider two practical problems of interest. First, we consider a device to device (D2D) communications setting, and the goal is to design a user trajectory while maximizing its connectivity to the internet. This problem is of vital interest, due to the surge in data-intensive applications in smartphones, and the consistent internet connectivity is becoming essential. For this problem, we consider a pair of pedestrians connected through a D2D link for data exchange applications such as file transfer, video calling, and online gaming, etc. The second problem is associated with the online planning of energy-efficient trajectories for unmanned surface vehicles (USV) under strong disturbances in ocean environments with both static and dynamic goal locations. We consider an unmanned surface vehicle (USV) operating in an ocean environment to traverse from start to destination. Different from the state of the art trajectory planning algorithms that entail planning and re-planning the full trajectory using the forecast data at each time instant, the proposed algorithm is entirely online and relies mostly on the current ocean velocity measurements at the vehicle locations. The detailed simulation results demonstrate the significance of the proposed algorithm on synthetic and real data sets. Video on the real world datasets can be found at \colr{https://www.youtube.com/watch?v=FcRqqWtpf\_0}}  

\end{abstract}

}
\maketitle

\IEEEdisplaynontitleabstractindextext
\IEEEpeerreviewmaketitle


\section{Introduction}
Trajectory design for motion planning is a core requirement for all autonomous systems. The foremost goal in the motion planning problems is to determine the optimal trajectory starting and ending at specified points, while also satisfying application-specific constraints, such as (i) avoiding clutter and other obstacles  \cite{coll-avoid}; (ii) maintaining connectivity \cite{li2018joint};  (iii) adhering to vehicle-specific constraints like turn radii or spatial envelop of the vehicle \cite{vhcap}; and (iv) having high energy efficiency \cite{eesto}. However, in practice, various environmental variables are time-varying and unknown in advance, motivating the need for online trajectory design. Herein, at any given time, the agent makes local observations, and uses only current available information to determine the next feasible waypoint. Moreover, the agent does not rely on predictions or forecasts, and the resulting trajectory is robust to arbitrary changes in the environment. In order to further motivate the online setting, we discuss two contemporary applications: rate optimal trajectory planning for device-to-device (D2D) communications in cellular networks and energy-efficient trajectory design in an ocean environment. 

\begin{figure}
	\subfigure[]{\includegraphics[width=0.46\linewidth, height = 0.31\linewidth, trim={0cm 0cm 0cm 0cm},clip]{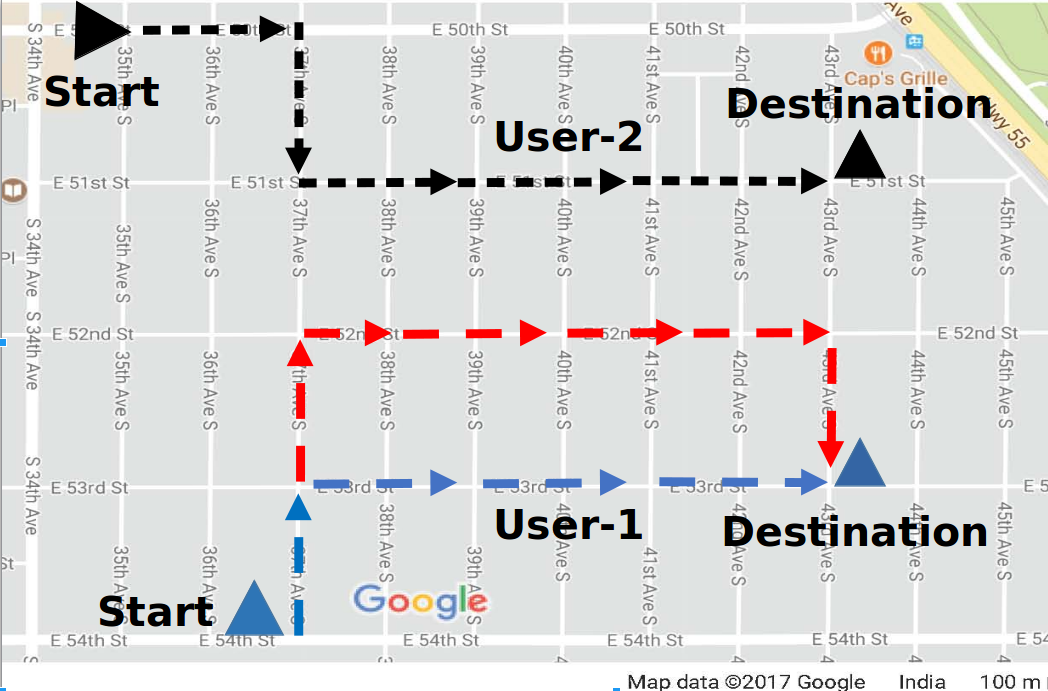}} 
	\subfigure[]{\includegraphics[width=0.54\linewidth,height = 0.35\linewidth, trim={0cm 0cm 0cm 0cm},clip]{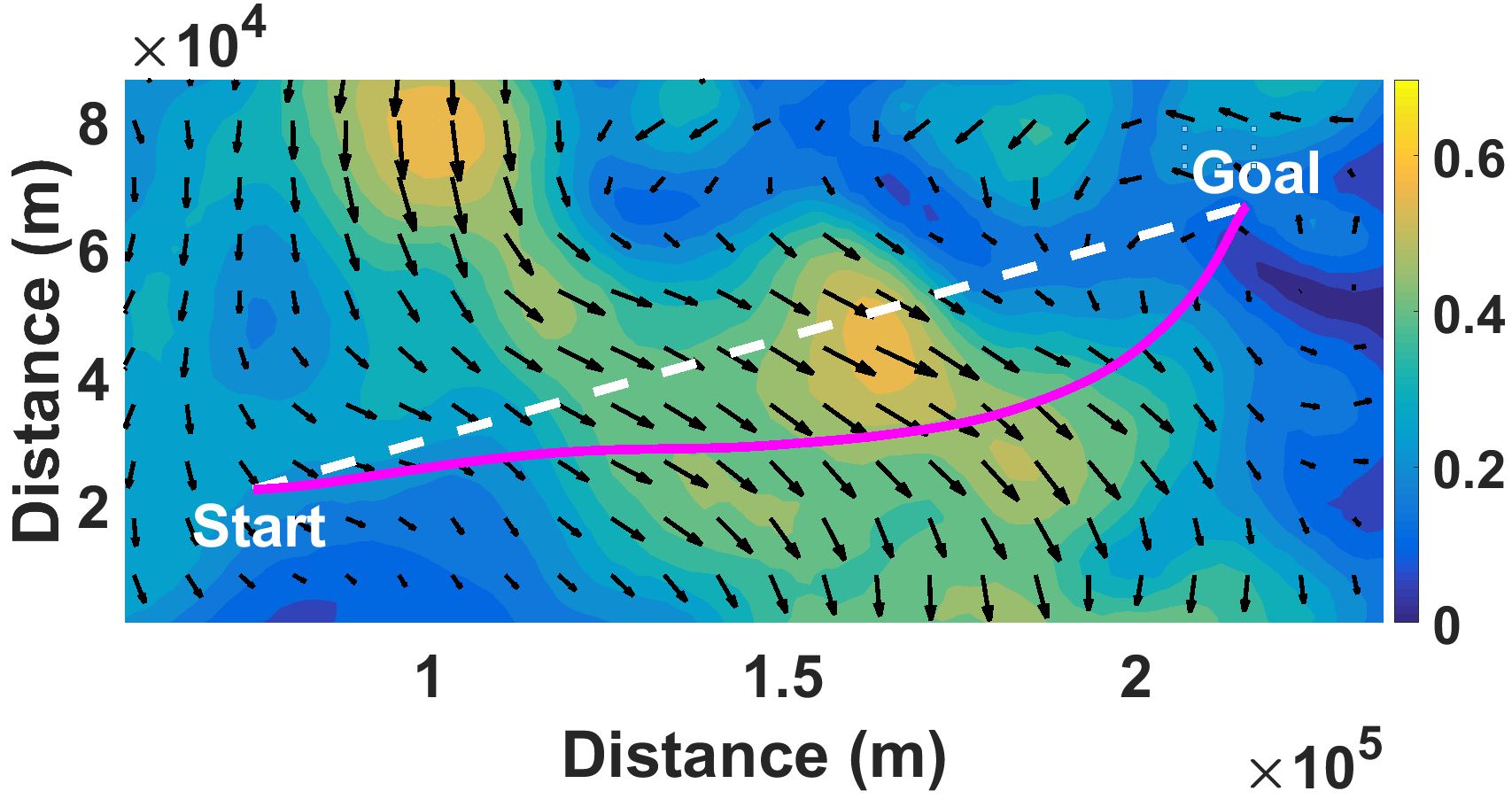}}\\
	\captionsetup{font=scriptsize}
	\caption{(a) A simple fixed hotspot model, blue dashed line shows a direct path while red dashed line shows an alternative path to increase data rate. (b) An instance of the ocean environment: the white dashed line shows the direct path, whereas the solid magenta line represents the energy-efficient path which utilizes the information of ocean current.}   
	\label{system_model} 
	\vspace{-0.2cm}
\end{figure}

\noindent \emph{\textbf{Rate-optimal trajectories:}}  The human need to be constantly connected has fueled the growth of social media platforms and associated smartphone applications.    
With the advent of data-intensive applications, the service providers now face the challenge of providing seamless connectivity `on the go' \cite{li2009future,chin2014emerging}. The need to maintain high spectral efficiency has prompted the researchers to look beyond the traditional cellular architecture and develop innovative features such as device-to-device (D2D) communications \cite{lin2000multihop}. Enabling direct communications between nearby cellular users has the potential to not only improve spectrum utilization, and throughput, but also enable disruptive peer-to-peer applications and services \cite{ji2016wireless}. Indeed, the D2D paradigm allows users to sustain connectivity at low costs and without overloading the cellular network \cite{asadi2014survey}. Thanks to the content heavy nature of the modern social media platforms, the cybercitizens of today are increasingly willing to modify their behavior in order to stay connected. For instance, urban areas have seen an unprecedented increase in the number of wifi hotspots \cite{mota2013feasibility}. As more D2D devices surface, it is likely that the users will be willing to modify their daily commute to stay connected. For instance, pedestrians may be willing to increase their commute times by a fraction so as to not only exchange more data and connect longer, but also sustain a higher quality of service. Fig. \ref{system_model}(a) shows the direct (blue) and  the alternative (red) path taken by the user in order to sustain higher data rate with a nearby mobile user (black).  \\
\noindent \emph{\textbf{Energy-efficient trajectories:}} Energy efficiency is a central issue in naval and aerial path planning tasks where the surrounding air/water flow may allow the vehicle to conserve energy. Designing energy-efficient trajectories is important, for instance, in ocean environments where the vehicles are often deployed for long-term autonomous missions such as surveying, mine hunting, chasing seaborne targets, oceanographic research, etc., with limited energy resources on board \cite{longterm}. Classical optimal trajectory design approaches in this context require current ocean information between the starting point and the goal \cite{app}. Here in Fig.\ref{system_model}(b) we display an instance of ocean environment where we show the direction of ocean currents and the longer trajectory (magenta) which incurs significantly lower energy as opposed to the straight line path. However, unless the trajectories are short, such information is only available in the form of low-resolution historical data or forecasts. For instance, approximate forecasts available from Regional Ocean Modeling System (ROMS) \cite{roms} or Navy Coastal Ocean Model \cite{ocean.model} provide ocean current velocities at points separated in the order of kilometers. Further, since the forecasts are obtained from satellite or field observations, they are available either at very coarse resolutions or only for specific geographical regions. Owing to the uncertainty in the ocean current forecasts, offline trajectory design of autonomous vehicles is potentially dangerous due to the possibility of being swept into shipping lanes or land \cite{uncert.ocean}. To counter these challenges, state-of-the-art approaches rely on re-planning trajectories using measured ocean current values \cite{eesto}. The optimal trajectory is first found using the forecast data, and subsequently at each time slot, forecasts are regenerated using measured velocities and the trajectory designs are revised. It is remarked that although such algorithms are capable of handling time-varying environments and a moving goal, the reliance on forecasts and high computational costs make them impractical for many systems.

This paper formulates a general utility-optimal trajectory design problem in the unknown and time-varying environments and develops an online algorithm for the same. Following the spirit of online convex optimization (OCO), the proposed algorithm, given the start location, generates the subsequent waypoint locations in an online fashion using an (projected) inexact online gradient ascent (IOGA) algorithm. Since each update step only uses the information at the current time instant, both the environment as well as the goal location are allowed to be time-varying. The performance of the IOGA algorithm is studied for the time-varying problem at hand and the resulting regret against the offline solution is shown to be sublinear. In order to encourage utility-optimal paths, the agent is allocated additional delay budget and the ensuing trade-off is explored via simulations. 
 
The applicability and generality of the proposed framework is demonstrated via two complementary applications, namely data rate-maximizing trajectories for commuting users and energy-efficient trajectories in ocean environments. In the former case, the goal is to maximize a user-specific utility function depending on the download data rate or on the signal strength while ensuring that the users reach their destinations. Specifically, we consider a pair of mobile users following a specific trajectory, say from a starting point to a destination. The users exchange data on a D2D link and need to sustain high data rates. The users are willing to take detour from the most direct path from the starting point to their respective destinations. The goal is to design an optimal trajectory for the users to not only maximize their data rates but also reach their respective destinations within the allocated excess delay budget. In the latter case, the goal is to design an energy-efficient watercraft trajectory that is capable of conserving energy by utilizing the ocean currents. At each time, only the subsequent waypoint is determined, while ensuring that certain maximum velocity and acceleration constraints are satisfied. An excess delay budget is also specified that again allows the agent to conserve additional energy \cite{timenr,witt}. We remark that the design of the utility functions and the adaptive step-size rules for each of the two algorithms is not straightforward and form the main contribution of the application sections.  

\noindent\textbf{Contributions:} \colb{The main contributions of this work include (a) development and analysis of the generic time-varying utility-optimal trajectory design problem; (b) performance analysis of the proposed IOGA algorithm for the case of time-varying objective and constraint functions, yielding a sublinear regret; (c) formulation and detailed study of the online trajectory design problem for a watercraft operating under strong ocean currents; (d) solving D2D trajectory planning problem using IOGA, by considering the uncertainty in user's positional estimates. This problem is inspired from \cite{amrit}, where the D2D trajectory optimization problem is formulated for the first time and solved using online gradient descent (OGD) under noiseless feedback ignoring real world uncertainties. The problem formulation in this work is different from \cite{amrit} and is considered within a more general context where more realistic time-varying constraints are now included. In addition, we propose to solve the problem in an online fashion under noisy gradient feedback and refer to it as IOGA. The modification required to bound the resulting regret of IOGA is not trivial and culminates in additional cumulative gradient and error variation terms in the expression for regret. Additionally, we derive lower bound to the offline dynamic regret that helps in deciding the optimality of any online policy. Next, the IOGA parameters are designed for the more sophisticated energy-efficient trajectory design problem not considered in \cite{amrit}. The performance of the proposed algorithm is tested on real-world ocean current data and is shown to be significantly better than the state-of-the-art algorithms. A special case of this work i.e., IOGA with noiseless gradient feedback (i.e., OGA) with application to planing energy efficient trajectories in time-varying ocean environments has appeared in \cite{icra}.}

The rest of the paper is organized as follows. Sec. II briefly reviews some related literature. In Sec. III we detail the general trajectory planning problem formulation besides assumptions and conditions under which desired regret bounds hold true. Sec. IV includes the derivations of regret bounds and analysis. Sec. V. details D2D and the energy efficient trajectory planning problems, their formulations and solution methodologies. Experimental validation is provided in Sec. VI. Finally, Sec. VII concludes the paper.

   
\section{Related Work}  
The problem of trajectory planning for a moving robot in a specified area has been widely studied in the robotics and control literature \cite{vasquez2016novel,zeng2015survey,mac2016heuristic}. For the time-varying setting that is of interest here, the standard approaches are based upon graph theory \cite{takei2010practical,likhachev2009planning,ferguson2006using}, Markov decision process based techniques \cite{fakoor2016humanoid,spaan2004point}, and model predictive control \cite{falcone2007linear,kim2014model}. Other advances in the path planning and robotics areas include dealing with uncertainties, dynamics, or multiple robots \cite{van2011lqg,liu2011coordinated,plaku2010motion}. The standard approach here is to generally formulate a complicated optimization problem and solve it using interior point method at each time step $t$. \colb{In contrast, the present paper specifically targets the problem from an online learning perspective, leading to a significantly simpler algorithm that is also amenable to well-defined performance guarantees without ignoring real world uncertainties.} In a similar vein, the path planning problem has also been considered under the aegis of reinforcement learning. Within this context, the problem is formulated as that of finding a sequence of feasible actions that take a robot from a source to destination \cite{al2013wind,konar2013deterministic}. Different from these works, { we do not assume that the distance to the next position from the current position is known to the moving user.} The problem of trajectory design has also been considered from a variational perspective, with the trajectory given by a continuous time function. The variational problems are often solved numerically via discretization and consequently are not amenable to the regret bounds such as those developed here \cite{fraichard1998trajectory}.

Coming to energy-efficient planning in ocean environments, graph-based search methods are designed to leverage the coarse resolution of data available 
\cite{headingangle,huynh2015,p47}. However, traditional approaches suffer from two general problems: (i) the quality and complexity of the solution is controlled by the selected discretization and (ii) the predicted trajectories often result in increased control costs because of the need to make sharp turns. Sampling-based methods represent the trajectory as a sample from a Gaussian process \cite{gpmp}, or sample
a set of noisy trajectories about the initialized trajectory \cite{stomp}, \cite{eesto}. Of these, \cite{eesto} also samples in time, thereby providing additional flexibility in the design of trajectories. An exploration parameter is often required to control the variance of the sampled trajectories, and a higher exploration parameter often leads to better but possibly more skewed trajectories. On the flip side, such sampling-based algorithms are highly sensitive and output a different trajectory for every run of the algorithm. Approaches such as \cite{eesto} and \cite{stomp} also require dense sampling rendering them slow in the presence of multiple constraints. Overall, the per-iteration run-time of sampling-based algorithms is at least $\mathcal{O}(KN^3)$ where $K$ is the number of samples, and $N$ is the total number of waypoints. In contrast, the proposed algorithm has only a single tuning parameter, no random components, and a run-time of $\mathcal{O}(N)$.

Optimization methods approach the trajectory design problem more directly. Offline approaches include methods utilizing parallel swarm search \cite{witt}, level-set expansion methods \cite{subramani,lolla}, and gradient-based approaches \cite{kruger2007optimal}. These approaches rely on forecast data and cannot be readily adapted to online and time-varying settings. Online approaches to trajectory design include the sequential convex optimization framework in \cite{trajopt} and covariant gradient ascent algorithm \cite{chomp}. However, existing online approaches require static settings and very little analysis exists for the time-varying setting at hand. Of these, most works provide time-optimal trajectory designs \cite{kruger2007optimal}. In contrast, the travel taken of the sampling-based algorithm such as \cite{eesto} may be different every time depending on the sampled trajectories and ocean current velocities. Finally, two separate optimization problems, one for minimizing the travel time, and another for minimizing the energy, are posed in \cite{p47}. Although \cite{p47} is demonstrated under time-invariant flow fields, it can be easily extended to time-varying scenarios. However, extensions to time-varying flow fields limit the accuracy of flow model used to predict the future temporal variations. In contrast, the current work  designs trajectories in an online fashion without such models.

The current work formulates the trajectory design problem within the OCO framework, while allowing time-varying objective function and constraints. Online learning algorithms such as online gradient descent (OGD) have a rich history and are generally studied in terms of their regret performance, where the performance of the algorithm is compared with that of a static benchmark \cite{UnconOCO,UnconOCOimp}. For dynamic environments however, it makes sense to consider dynamic regret, where the benchmark itself is time-varying \cite{UconOCOdyn1,UconOCOdyn2}. The desiderata here is that over a horizon of $T$ time slots, the regret should be a sublinear function of $T$ if the cumulative drift of the benchmark is also sublinear. \colb{More importantly, recent studies have shown interest in analyzing inexact gradient like methods \cite{bedi2017adversarial,8598992,8268092} as calculating approximate gradient is often cheap in machine learning literature. }

\colb{To this end, in the present work, we analyze the inexact online gradient ascent algorithm considering the so called offline regret, where the benchmark is the optimal trajectory designed with full knowledge of the future \cite{regret1}.} Similar regret metrics have recently been considered for problems with separable constraint functions \cite{adaptiveOL,OCOtiv3, neely,adapoco1}. Different from these work, the problem considered here comprises of generic temporally coupled constraint functions. \colb{To the best of our knowledge, this is the first work to consider offline benchmark under time-varying coupled constraints.} 
\section{General Trajectory Optimization Problem} \label{GTOP}
\colb{This section formulates the general utility-optimal trajectory planning problem and puts forth an inexact online gradient ascent algorithm (IOGA) for the same. The term \textit{inexact} is used here, as the gradient feedback is noisy.} We begin with formulating a constrained trajectory utility maximization problem that would serve as a benchmark for the proposed trajectory planning algorithm. Subsequently, the online formalism is used to develop the proposed algorithm and establish that it achieves sublinear regret. The algorithm and formulation in this section build upon the trajectory design problem considered by \cite{amrit} in the context of communication networks. The analysis here considers a more general class of time-varying constraints and is therefore applicable to a wider variety of trajectory planning problems. \colb{For instance, consider the energy-efficient trajectory planning problem in time-varying ocean environments, a constant velocity assumption can be relaxed to generate more energy-efficient trajectories, with the help of a time-varying constraint formalism (see \eqref{constr}).}  
\subsection{Trajectory Utility Maximization}
The trajectory planning problem is expressed as that of finding the optimal agent locations $\x(t)$ at each $t\in\N$ subject to various restrictions. In the literature, these locations are also referred to as waypoints that an agent must traverse and the trajectory is simply the collection $\{\x(t)\}_{t\geq1}$. Here, the time is measured in units of `slots' whose lengths may otherwise be arbitrary. At time $t = 1$, the current location of the agent $\x(1) = \mathbf{s}$, as well as the goal (destination) location $\d(t)$ are given. The framework subsumes the multi-agent scenario where $\x(t)$ is the stacked version of the individual agent locations $\{\x_i(t)\}_{i=1}^M$ at time $t \in \N$. Let $T_{\text{ETA}}$ denote the expected time-of-arrival (ETA) for the agent when it follows the straight line path between $\x(1)$ and $\d(1)$ calculated either using forecast data or the prevailing environmental conditions at time $t = 1$. However, due to the dynamic nature of the problem, at time $t = T_{\text{ETA}}$, (a) the agent need not be at $\d(1)$ if it follows the straight line path and (b) the goal location would have already changed to $\d(T_{\text{ETA}})$. To this end, a budget of $\delta$ excess time slots is available and may be used to traverse a more utility-optimal trajectory. The general utility-optimal trajectory design problem may be formulated as follows 
\begin{subequations}\label{trajopt}
\begin{align} 
\left\{\x^\star(t)\right\}_{t=1}^T &=\arg\!\!\!\!\!\! \max_{\left\{\x(t)\right\}_{t=1}^T} \  \sum_{t=1}^{T} U_t(\x(t)) \label{trajopt1}\\
 \text{s. t. }\ & g_t(\x(t),\x(t+1)) \leq 0 \hspace{0.5cm}  1\leq t \leq T-1 \label{trajopt2} \\
 & \x(1) = \mathbf{s} \\ 
 & \x(t) \in \mathcal{X} \hspace{2.37cm}  1\leq t \leq T  \label{trajopt3}
\end{align}
\end{subequations}
where $T = T_{\text{ETA}}+\delta$, $\X\subset\R^{2}$ is non-empty closed and convex set that represents the {Euclidean space} in which the agent moves, $U_t:\X \rightarrow \R$ is a time-varying concave utility function, and $g_t:\X \times \X \rightarrow \R$ represents the convex coupling constraint between successive way-points. Both $U_t$ and $g_t$ may depend on other physical parameters of the agent/environment such as the goal location $\d(t)$, maximum agent speed $v_{\max}$ (measured in meters per time-slot), velocity of the ocean currents etc. While the subsequent analysis will require that the goal trajectory $\{\d(t)\}_{t=1}^T$ and the functions $U_t$ and $g_t$ to be slowly time-varying, their temporal variations are otherwise arbitrary and possibly adversarial. In the multi-agent scenario, we have that $\X \subset \R^{2M}$ and $g_t$ is replaced with the vector valued function $\g_t:\X \times \X \rightarrow \R^M$. For the sake of brevity, the subsequent analysis will be developed for a single agent scenario. However the online algorithm and the relevant regret bounds can be readily obtained for $M \geq 1$. 

Problem \eqref{trajopt} may be solved in an offline fashion using any existing convex optimization algorithm. In practice however, decisions at time $t$ can only be based on information available till time $t$. In the spirit of online convex optimization (OCO), we adopt the adversarial learning model, where the learning process can be viewed as a game between the agent and the environment. \colb{ The agent located at $\x(t)$ at time $t$ obtains $(\gradt U_t, g_t)$ from the environment and makes the decision to move to a new location $\x(t+1)$. Observe that the inexact gradient of the utility function $\gradt U_{t+1}(\x(t+1))$ is revealed only after the agent moves to $\x(t+1)$. The online framework requires no assumptions regarding the stationarity of environmental effects. Instead, online algorithms are always characterized with regards to their worst-case performance.} 

\subsection{Trajectory Optimal IOGA}
Online gradient-based algorithms have been widely studied for unconstrained optimization problems. For the temporally coupled problem at hand, we adopt the policy of online gradient ascent (OGA) that uses bounded stochastic gradient feedback. Altogether, we consider the case where the agent doesn't have the full functional form of $U_t$, instead a noisy gradient feedback is available. This is generally the case in practical scenarios; for instance, while planning energy-efficient trajectories in time-varying ocean environments, the utility function becomes noisy due to the uncertainty in the ocean current measurements.

Note that the application of the OGA algorithm is tricky, since a step in the direction of the gradient may render the iterate infeasible. To this end, we put forth the variable step-size IOGA algorithm that utilizes the bounded stochastic gradient feedback and, takes the form
\begin{align}\label{online}
\xh(t+1) &= \px\left(\xh(t) + \tfrac{1}{\gt}\gradt U_t(\xh(t))\right)  
\end{align}
for $t = 1, 2, \ldots, T-1$, where $\gt$ is the step-size or the learning rate, $\gradt U_t(\cdot)$ is the inexact gradient i.e., $\gradt U_t(\cdot) = \nabla U_t(\cdot) + \nt$ for some noise $\nt\in\R^2$, and the initialization $\xh(1) = \mathbf{s}$. Different from the classical diminishing or constant step-size online gradient algorithms, the key idea here is to choose $\gt$ for each $t\geq 1$ such that the iterates stay feasible to satisfy \eqref{trajopt2}. For instance, consider the case when $\X = \R^2$ and  $g_t(\x,\x') = \norm{\x-\x'}- v_{\max}$ for any $\x$, $\x' \in \R^2$. Then, the gradient boundedness property $\norm{\gradt U_t(\x)} \leq \tilde{G}$ implies that 
\begin{align}
\norm{\xh(t+1)-\xh(t)} \leq \frac{\tilde{G}}{\gt} \leq v_{\max}
\end{align}
if we choose $\gt \geq \tilde{G}/v_{\max}$ for all $1\leq t\leq T$. \colb{It is remarked that the proposed algorithm and the analysis is only applicable to scenarios where such a choice of $\gt$ exists. In other words, the iterate feasibility condition must be explicitly verified for the problem at hand. Also, note that the proposed IOGA algorithm subsumes OGA when the gradient feedback is noiseless. Therefore, IOGA is more general and hence pertinent to realistic scenarios as opposed to OGA. To this end, the agent is expected to make use of the inexact gradient to learn and improve the trajectory estimates while adhering to the time-varying constraints.} Before proceeding with the performance analysis, a remark regarding function availability is due. 

\begin{rem}
Within the online setting considered here, only the gradient to the current utility function $U_t$ is revealed at time $t$. In many applications however, reliable predictions of the environment at time $t+1$ may be available at time $t$. For such settings, it may be better to use the modified updates that take the form
 \begin{align}\label{online2}
 \xh(t+1) &= \px\left(\xh(t) + \tfrac{1}{\gt} \gradt U_{t+1}(\xh(t))\right)
 \end{align}
for $t = 1, 2, \ldots, T-1$. It will be subsequently  established that the updates in \eqref{online2} performs better than the one in \eqref{online}.
\end{rem}

\subsection{Performance Analysis}
Having detailed the proposed trajectory optimal OGA and IOGA algorithms, we are now ready to analyze each of its performance. We begin with stating the necessary assumptions. 
\begin{itemize}
	\item [({A1})] \textbf{Strong concavity:} the utility function $U_t$ is $\mu$-concave or equivalently, $U_t(\x)+\frac{\mu}{2}\norm{\x}^2$ is concave.  
	\item [({A2})] \textbf{Lipschitz continuous gradient:} the utility function $U_t$ is $L$-smooth, i.e., $\norm{\nabla U_t(\x)-\nabla U_t(\x')}\leq L \norm{\x-\x'}$ for all $\x$, $\x'\in\X$.   
	\item [({A3})] \textbf{Bounded gradients:} the utility function gradients are bounded, i.e.,  $\norm{\nabla U_t(\x)}\leq G$ for all $\x\in\X$. 
	\item [({A4})] \textbf{Feasibility:} the IOGA iterates $\xh(t)$ adhere to \eqref{trajopt2}.     
	\item [({A5})] \colb{ \textbf{Error bound: }  bounded second order moments i.e., $\E[\norm{\nt}^2] \leq \varepsilon_t^2 < \infty$}.                     
\end{itemize}
\colb{Observe here that while Assumptions (A1)-(A3) and (A5) are standard and are widely used, Assumption (A4) is specific to the problem at hand. As stated earlier, (A4) must be explicitly verified and is generally met by choosing $\gt$ appropriately. We emphasize that we do not require the stochastic error to be stationary or unbiased.} 

These set of assumptions enable us to characterize the performance of the proposed trajectory optimal IOGA algorithm through the notion of regret, which compares the utility obtained by the proposed algorithms against that obtained by an adversary. Note that the adversary has complete information about the future and solves \eqref{trajopt} in an offline manner. An algorithm is termed no-regret, if the regret grows sublinearly with $T$. 

For the problem at hand, we consider the so-termed \emph{offline dynamic regret} \cite{regret1}, that can be written as
\begin{align}\label{regret}
\textbf{Reg}_T:= \underbrace{\left[\s U_t(\x^\star(t))\right]}_{\text{offline}} - \underbrace{\left[\s U_t(\xh(t))\right]}_{\text{online}}.
\end{align}
Observe that the definition in \eqref{regret} is different from the classical static regret where the benchmark is a static $\x^\star$ that solves \eqref{trajopt1}. Observe that, static benchmark has no meaning here as it would suggest to stay at the same location for all the time. Instead, the offline regret contains a different $\x^\star(t)$ for each $U_t$ and is therefore stronger. Indeed, it is known from \cite{regret1} that for optimization problems with temporally coupled constraints, the offline regret is generally linear in $T$. 

\subsection{Parameter Variations and Error}
The regret bound defined in \eqref{regret} is calculated in terms of various parameter variations pertaining to different cases analyzed. We propose to study the performance of the IOGA algorithm taking the form of updates as defined in \eqref{online}. We now, define those variations, starting with the squared path length of the adversary, defined as 
\begin{align}\label{st}
\St :=\sum_{t=1}^{T-1} \norm{\x^\star(t+1)-\x^\star(t)}^2.
\end{align}
The squared cumulative gradient variation as 
\begin{equation}
\Gt := \sum_{t=1}^{T-1}\underset{\x\in\X}{\max}\norm{\nabla U_{t+1}(\x) - \nabla U_t(\x)}_2^2.
\end{equation}
Intuitively, in the former case large variations in $\x^\star(t)$ make it difficult for the agent to follow, and may lead to a linear regret. For instance, if the squared path length is linear in $T$, e.g., if the goal location $\d(t)$ moves away from the agent by $v_{\max}$ meters per time slot, the agent may never be able to catch up with the goal. The cumulative gradient variation is also similar and has been widely used in the context of online learning and dynamic optimization \cite{cgv1,cgv2}. \colb{Next, the cumulative squared error bound is defined as,
	\begin{align}
	E_T &= \sum_{t=1}^{T} \varepsilon_t^2. 
	\end{align}  
	As we shall show later, the regret will be bounded by a function of the cumulative squared error. Such a behavior is also expected, since arbitrary gradient errors can accumulate over iterations. However, it will be established that regardless of the error statistics, their impact cannot be catastrophic; small gradient errors will still have a small impact on the regret. }

\section{Regret bounds and analysis}
Let $\gmax:=\max_{1\leq t\leq T-1} \gt$ and $\gmin:=\min_{1\leq t \leq T-1}\gt$ be constants that do not depend on $T$. The following theorem summarizes the regret bound for \eqref{online}. 

\begin{thm}\label{thm}
	Under the assumptions {(A1)-(A4)} and for $\gmax < 2\gmin - L$, the sequence of $\xh(t)$ generated by IOGA update in \eqref{online} adheres to the expected regret bound 
	\begin{align}\label{regioga}
 \colb{\E[\textbf{Reg}_T] = \O{\sqrt{T(\St + \Gt + E_T)}}.}
	\end{align}
\end{thm}



The result in Theorem \ref{thm} states that for large values of $T$ and for a sub-linearly time-varying adversary, the trajectory optimal IOGA algorithm incurs a sub-linear regret. As a simple example, consider the case when the optimal path is the straight line from $\mathbf{s}$ to $\d(T)$ covered in $T$ time slots. Since a distance of $\norm{\d(T)-\mathbf{s}}/T$ is covered at every time slot, the squared path length $\St$ is $\O{1/T}$ if the target remains within a fixed bounded area. Assuming $\Gt$ and/or $E_T$  to be of the same order, it can be seen that the regret bound shown in \eqref{regioga} is constant for this example. In other words, the IOGA algorithm for this example is only suboptimal by a constant regardless of how high $T$ is. Equivalently, the per-time instant utility loss, arising from not knowing the future, tends to zero as $T$ increases. Regret will generally be higher if the optimal path is not a straight line.  

\colb{For the regret bound to be meaningful, it is necessary that $E_T$ be sublinear in $T$, a common requirement when studying online algorithms \cite{bedi2017adversarial, 8598992}. Sublinearly accumulating error is also natural, for instance, in the ocean environment, where noise filtering techniques can be used to drive the systematic measurement errors to zero, at least after a few iterations. Further, the dependence on $E_T$ instead of the individual terms $\varepsilon_t$, allows for intermittent sensor failures leading to spikes in the measurements. As long as the such spikes occur rarely so that $E_T$ remains sublinear and non-dominant, they would not effect the overall regret rate.}

\begin{IEEEproof}[Proof sketch]
The proof of Theorem 1 follows along the lines of that in \cite{amrit} but includes modifications required to handle the generic time-varying convex constraint function $g_t$ and the noisy gradient feedback. It is remarked that the modification from \cite{amrit} is not trivial and changes the proof as well as the final result considerably. \colb{Also note that the condition $\gmax < 2\gmin - L$  
	implicitly implies that $\gmin > L$ and hence it holds that $\gmax > L$. The details of the complete proof are provided in section I of the supplementary material. The main idea of the proof is to utilize the the quadratic lower and upper bounds, thanks to smoothness and strong concavity of $U$, in order to establish a fundamental inequality. This inequality is subsequently used along with assumptions (A3)-(A5) and the Peter-Paul inequality to obtain the required result.} 
\end{IEEEproof}
 
The result in Theorem \ref{thm} is stronger than similar results that have been shown to hold  for dynamic regret in \cite{mokhtari2016online}. Likewise, the offline regret is also known to be linear in general \cite{regret1}. In contrast, Theorem 1 establishes a generic sublinear bound on the offline regret under relatively mild assumptions. 

Also, the regret for the modified IOGA updates as in $\eqref{online2}$, takes the form 
\begin{align}\label{reg4}
\colb{\E[\textbf{Reg}_T] = \O{\sqrt{T(\St+E_T)}}.}
\end{align}
The proof for \eqref{reg4} follows along the same lines as that of \eqref{regioga}. Specifically, the proof begins in the same way but does not require the introduction of the term $G_t(\xh(t))$, finally resulting in the regret expression that does not contain $\Gt$. The detailed proof is provided in section II of the supplementary material attached to this manuscript. 

\colb{We remark here that the IOGA subsumes the noiseless gradient feedback case and we refer to it as OGA. The corresponding regret bounds for the updates shown in \eqref{online} and \eqref{online2} are given by
\begin{align}
\textbf{Reg}_T & = \O{\sqrt{T(\St + \Gt)}},   \label{reg3} \\
\text{and} \hspace{1cm}\textbf{Reg}_T & = \O{\sqrt{T(\St)}}, 
\end{align}
respectively, obtained by setting $E_T = 0$.} The subsequent remarks discuss the choice of problem parameters and the optimality of the bounds in Theorem \ref{thm}. 

\begin{rem}
\colb{Assumption (A1) is integral to obtaining a sublinear offline regret bound. The present analysis does not apply to concave utility functions since $\mu > 0$ is utilized to balance the error due to parametric variations; for further details see proof of Theorem \ref{thm} from section I of the supplementary file. Indeed, for the present proof, we require that
\begin{equation}
\mu > \frac{1}{\sqrt{T}}\sqrt{o\left(\max\{\St, \mathcal{G}_T,E_T \} \right)},
\end{equation}   
to yield a sublinear regret. It remains open to see if (A1) can be relaxed to allow weaker conditions like error bound \cite{bedi2017adversarial}}. 
\end{rem}
 
\begin{rem}
\colb{The offline regret defined in \eqref{regret} cannot be sublinear unless the squared path length is sublinear. To see this, it suffices to establish that there exists a sequence of functions $\{(U_t,g_t)\}_{t=1}^T$ with linear $\St$ for which no policy can have a sublinear regret.} 
	
\colb{To construct such a sequence of functions, we consider the scalar case for simplicity. Let $\{w(t)\}_{t=1}^T$ be some parameters such that $w(t) \in [-W,W]$ is revealed by the adversary after the action $x(t)$ has been taken. We consider the class of functions $U_t(x(t)) = -\frac{1}{2}(x(t)-w(t))^2$ and $g_t(x,x') = \abs{x-x'}^2-4W^2$. For this choice, it can be seen that if $x(t) \in [-W,W]$, the constraint is always satisfied. In fact, the optimal solution is given by $x^\star(t) = w(t)$ and hence $\Gt$ and $\St$ can be bounded as:
\begin{align}
\St= \Gt = \sum_{t=1}^T(w(t+1)-w(t))^2 \leq 4W^2T.
\end{align} 
The regret is given by 
\begin{align}
\textbf{Reg}_T = \frac{1}{2}\sum_{t=1}^T (x(t)-w(t))^2.
\end{align}
Suppose that the adversary chooses $w(t) = -W\text{sign}(x(t))$. In that case, it can be seen that the regret lower bounded as
\begin{align}
\textbf{Reg}_T = \frac{1}{2}\sum_{t=1}^T (x^2(t) + W^2 + 2W\abs{x(t)}) \geq \frac{W^2}{2}T,
\end{align}
irrespective of the policy $x(t)$ used. Therefore, for a constant $W$ (e.g. $W = 1$), it is possible to have functions $\{(U_t,g_t)\}_t$ for which the regret is linear if no other restrictions are imposed on the squared path length. On the other hand, it can also be seen that for $W = T^{-\gamma}$ for $\gamma > 0$, the regret is at least $\O{T^{1-2\gamma}}$, which is lower than the upper bound of $\O{T^{1-\gamma}}$ obtainable using OGA. Such a gap suggests that there may exist online algorithms which can improve upon the regret rate in Theorem \ref{thm}.}
\end{rem}



\section{D2D trajectory planning problem}
This section considers the D2D trajectory optimization problem and demonstrates the applicability of the proposed algorithm for the same. 
 
\subsection{General Problem Formulation}
Consider a pair of users (say $x$ and $y$) located on the $\mathbb{R}^2$ plane and connected via a D2D link. Let the locations of the two users at discrete time $t \in \mathbb{N}$ be denoted by $\x(t),\y(t) \in \mathcal{S} \subset \mathbb{R}^2$, where $\mathcal{S}$ is the set of viable user locations. The origin and destination of the $i$-th user are denoted by $\mathbf{s}_i$ and $\d_i$ respectively, where $i \in \left\{x,y\right\}$. For the base case, the users do not modify their behavior and travel from $\mathbf{s}_i$ to $\d_i$ along the shortest path. When there are no obstacles, the shortest path is simply the straight line joining $\mathbf{s}_i$ and $\d_i$ and has length $\norm{\mathbf{s}_i-\d_i}_2$. On the other hand, when the user is only allowed to move along a grid, the length of the shortest path is the city-block or Manhattan distance $\norm{\mathbf{s}_i-\d_i}_1$. The distance between the source and the destination is henceforth denoted by $\norm{\mathbf{s}_i-\d_i}$, where the norm could be Euclidean, Manhattan, or any other convex distance metric. Also each user travels at the maximum speed of $v_{\max}^i$ units per time slot and therefore takes time $T_i^{ETA}:=\norm{\mathbf{s}_i-\d_i}/v_{\max}^i$ to reach the destination via the shortest path.

Here we consider a basic scenario where we plan the trajectory for a single user (say $\x$) such that the utility derived on the D2D link (say with user $\y$) will be maximum. However, the algorithm applies to the multi-user scenario as well where optimal trajectories for both the users are found. Recall  that the trajectory is a set of waypoints $\left\{\x(t)\right\}_{t\geq 1}$. More importantly, we assume a simple scenario where user $\y$ is cooperative and hence may wish to keep its current and future locations public. To this end consider a pair of users $(\x,\y)$ communicating on the D2D link and derive a utility of $U(\norm{\x(t)-\y(t)}_2)$ at time $t$. Here, $U$ is a non-increasing function of $\norm{\x(t)-\y(t)}_2$. Examples of utility functions may include the average received signal strength modeled as
\begin{equation}
U_{\text{RSS}}(\norm{\x(t)-\y(t)}_2) = RSS =\frac{1}{\norm{\x(t)-\y(t)}_2^{\alpha_p}},
\end{equation}
where $\alpha_p$ is the path loss parameter, or functions thereof. For instance, the average signal-to-noise ratio for the additive white gaussian noise channel with noise power $\sigma^2$ given by
\begin{equation}
U_\text{SNR}(\norm{\x(t)-\y(t)}_2)= SNR = \frac{RSS}{\sigma^2+RSS},
\end{equation}
and the channel capacity for a channel with bandwidth $W$ given by 
\begin{equation}
U_C(\norm{\x(t)-\y(t)}_2)=W\log_2(1+SNR),
\end{equation}
are some common examples. It is emphasized that in practice, the exact D2D rate would likely depend on a number of other factors, such as the application used, overheads, shadowing, etc. Given the inherent uncertainty in estimating the D2D rate in an urban setting, it may not necessarily be prudent to choose a complicated utility function. We set $T :=\min\{T_1,T_2\}$ and assume that the users disconnect as soon as one of them reaches the destination. Therefore, the cumulative utility obtained by the two users is given by $\sum_{t = 1}^{T} U(\norm{\x(t)-\y(t)}_2)$. 

As already mentioned in Sec. \ref{GTOP}, here we consider a scenario where the users are willing to take a longer path to their destinations in order to achieve a higher utility. Specifically, let $\delta_i$ be the excess delay user $i$ is willing to incur. In other words, it is acceptable for user $i$ to reach its destination in time $T_i = T_i^{ETA}+\delta_i$ so as to achieve a higher cumulative utility. Note that in the current formulation, the excess delay $\delta_i$ is an exogenous variable set by the user prior to starting. The D2D trajectory optimization problem that solves for the trajectory of user $\x$ can  be written as
\begin{subequations}\label{main}
	\begin{align}
	&\hspace{-2cm}\max_{\x(t)} \hspace{-3mm}\sum_{t = 1}^{\min_i\{T_i\}} U(\norm{\x(t)-\y(t)}_2), \nonumber\\
	&\hspace{-1.5cm}\text{s. t. } \ \x(1) = \mathbf{s}_{\x}, \label{start_rel}\\
	&\hspace{-0.7cm}\x(T_x^{ETA}+\delta_x)=\d_x ,  \label{dest_rel}\\
	&\hspace{-0.72cm}\norm{\x(t+1)-\y(t)}\leq \alpha(t) v_{\max}^x ,\nonumber\\
	&\hspace{1.5cm}1\leq t < T_x^{ETA}+\delta_x,\label{vel} \\
	&\hspace{-0.68cm}\x(t) \in \mathcal{X}. \label{set}
	\end{align}
\end{subequations}
Here, the constraints in \eqref{start_rel} and \eqref{dest_rel} ensure that the user trajectory begins and end at the specified starting and destination points. The constraint in \eqref{vel} enforces the maximum velocity constraint under the appropriate norm when $\alpha(t) =1$. Note that $\alpha(t) \in (0,1]$ scales the maximum velocity of the user by taking into account the current state of the environment.   
Finally, the constraint in \eqref{set} ensures that the user trajectory stays within the viable region. In an urban setting, the set $\X$ may represent the public areas such as roads, streets, alleys, etc. Such a constraint can generally be represented as a union of $K$ convex polygons, where each polygon represents an unobstructed area, such as the length and width of a road. Formally, each polygon is described by a set of affine inequalities of the form $\A_k\x(t) \leq \mathbf{b}_k$ for $1\leq k \leq K$ \cite{blackmore2006probabilistic}. It is emphasized that $\X$ is a union, not an intersection, of such polygons, and is generally non-convex. However, we assume $\X$ to be convex since developing performance guarantees for non-convex $\X$ is significantly harder and will not be pursued here. To this end, we consider a special case of \eqref{main} with the following two assumptions (a) the utility function $U$ is concave in $\x(t)$ for all $t$; and (b) the set of viable locations $\X$ is convex. The first assumption is satisfied, for instance, if the function $U(\cdot)$ is concave and non-increasing, e.g., $U(x) = -x^2$. On the other hand, the second assumption is not generally satisfied specifically in the areas with obstacles but is required to develop meaningful theoretical guarantees.  

In the general D2D setting, however, the offline formulation in \eqref{main} has limited applicability. To begin with, the formulation encodes a static scenario, where users are co-operative, and their destinations remain unchanged. In practice, the target location of the user may change after the user has already started, rendering the solution to \eqref{main} useless. More importantly, the users may like to interact online but not necessarily be willing to reveal their future and final destinations due to privacy and security concerns. The a priori unavailability of the full problem information motivates the need for online algorithms that are capable of handling time-varying and uncertain parameters. To this end, consider the case when $\y(t)$ is an exogenous variable that is revealed at each time $t$. Likewise, we also allow the final destination of user $\x$ to be a time-varying quantity $\d_x(t)$. The resulting time-varying problem can no longer be expressed in the form of \eqref{main}. Such problems have traditionally been considered within the context of rolling horizon planning or model predictive control.


Next, we reformulate the D2D trajectory optimization problem such that it is amenable to planning in an online fashion. The idea is to utilize the proposed IOGA algorithm by properly designing the objective and constraint functions such that they capture the requirements of the D2D trajectory planning application. We assume that the users are non-cooperative so that the current location of user~$\y$ given by $\y(t)$ and the current destination of user~$\x$ denoted by $\d_x(t)$ are exogenous and only available to user~$\x$'s device at time $t$. \colb{Different from \cite{amrit}, we do not assume the exact location of user~$\y$ to be known. Instead, $\y(t)$ is only approximately known, and the resulting uncertainty in the location estimates translates to the error in the gradient feedback.} For the sake of simplicity, we consider the utility function to be $U(x) = -x^2$. Since we optimize only the trajectory of user~$\x$, we drop the subscripts from the variables and parameters corresponding to user~$\x$. Therefore the new parameters of user~$\x$ become $v_{\max}$, $T = T_{ETA} + \delta$, $\mathbf{s}$, $\d(t)$, (iii) The maximum velocity in each slot is upper bounded by $v_{\max}$ and hence $\alpha(t) = 1$. \colb{We remark here that the extension to arbitrary concave utilities is left as future work.}
\subsection{Design of the Utility Function $U_t$} We consider the utility function for D2D application as follows
\begin{align}
U_t(\x(t))&:=-\lambda(t)\left\| \x(t)\!-\!\y(t)\right\|^2-(1-\lambda(t))\left\| \x(t)\!-\!\d(t)\right\|^2 \nonumber\\
&= - \norm{\x(t)-\bell(t)}_2^2 \label{loss}
\end{align}
where we define the leading path $\bell(t) := \lambda(t) \y(t) + (1-\lambda(t))\d(t)$. The sequence $\{\lambda(t)\}$ is included in order to control the relative importance assigned to being near user~2 and reaching the destination. Specifically, it is required that $\lambda(t)$ is a decreasing sequence that goes from 1 to 0 as $t$ goes from $1$ to $T$. Such a choice of $\lambda(t)$ ensures that the user places increasingly higher importance to reaching the final destination. Note that since $\lambda(t) \in [0,1]$, the optimization problem in \eqref{online} is still convex. Note that, differently from \eqref{main}, the constraint to reach the final destination is relaxed but is incorporated within the objective function.

We emphasize here that the cost function need not be the least-squares loss function in \eqref{loss}. More generally, it may be reasonable to choose convex functions of the form $U_t(\x(t)) :=- h(\norm{\x-\bell(t)})$ for some scalar function $h$. For instance, the square loss function that results in the gradient being proportional to $\x(t)-\bell(t)$. At the start, when the user is far from $\bell(t)$, the user may take larger steps, while it may slow down as it gets closer to $\bell(t)$. Another interesting choice is motivated from the Huber function, and takes the form
\begin{align} \label{huber}
\!\!\!\!h(d) = \begin{cases}\! \frac{1}{2}d^2 & d \leq v_{\max}\\
v_{\max}(1\!-\!\mu)d \!+\! \frac{\mu}{2}d^2 \!-\!  \frac{(1-\mu^2)v_{\max}^2}{2} & d > v_{\max}
\end{cases}
\end{align}
for some parameter $\mu$. Observe that the Huber function is same as the squared distance as long as the distance is less than $v_{\max}$. However, whenever the distance is larger than $v_{\max}$, the penalty is a convex combination of the linear and squared error penalties. The constants are adjusted to ensure that the function $U_t = -h(\norm{\x(t)-\bell(t)})$ adheres to (\textbf{A2}). As compared to the squared loss function, the Huber loss puts a smaller penalty when the user is far from $\bell(t)$. 

For this choice of objective function, its gradient can be written as
\begin{align}
\nabla U_t(\x)&= \begin{cases} \bell(t)-\x & \hspace{-3.2cm} \norm{\x - \bell(t)} \leq v_{\max},\\
v_{\max}(1-\mu)\frac{\bell(t)-\x}{\norm{\bell(t)-\x}} + \mu(\bell(t)-\x), \ \ o \symbol{92} w
\end{cases}\\
\nabla U_t(\xh(t)) &=\mu(\bell(t)\!-\!\xh(t)) \!+\! (1\!-\!\mu)P_v(\bell(t)\!-\!\xh(t)),  
\end{align}
where the projection operation is defined as
\begin{align}
P_v(\w)  &:= \arg\min_{\v} \norm{\v-\w}^2
&\text{s. t.} \norm{\v} &\leq v_{\max}.
\end{align}
\subsection{Design of the Constraint Function $g_t$} Now we define the constraint function in this context as follows
\begin{equation}\label{D2D_constraint}
g_t(\x(t),\x(t+1)) := \norm{\x{(t+1)}-\x(t)} \leq v_{\max}.
\end{equation}
As already described such a constraint enforces to satisfy motion constraints i.e., velocity of the user to be upper bounded by its maximum velocity achievable. We remark here that the solution to \eqref{main} and \eqref{trajopt} will generally not be the same. In particular
since the constraint in \eqref{dest_rel} is relaxed, the trajectory obtained by \eqref{trajopt} with the aforementioned objective and constraint functions, might stop short of the destination, depending on the manner in which $\lambda(t)$ is decreased. It will however be shown later in this Sec. that $\norm{\x(T)-\d}$ is generally small for the choice $\lambda(t) = t/T$. 
\begin{figure}
          \centering   
          \includegraphics[width=0.9\linewidth,height = 0.75\linewidth,trim={0cm 0cm 0cm 0cm},clip]{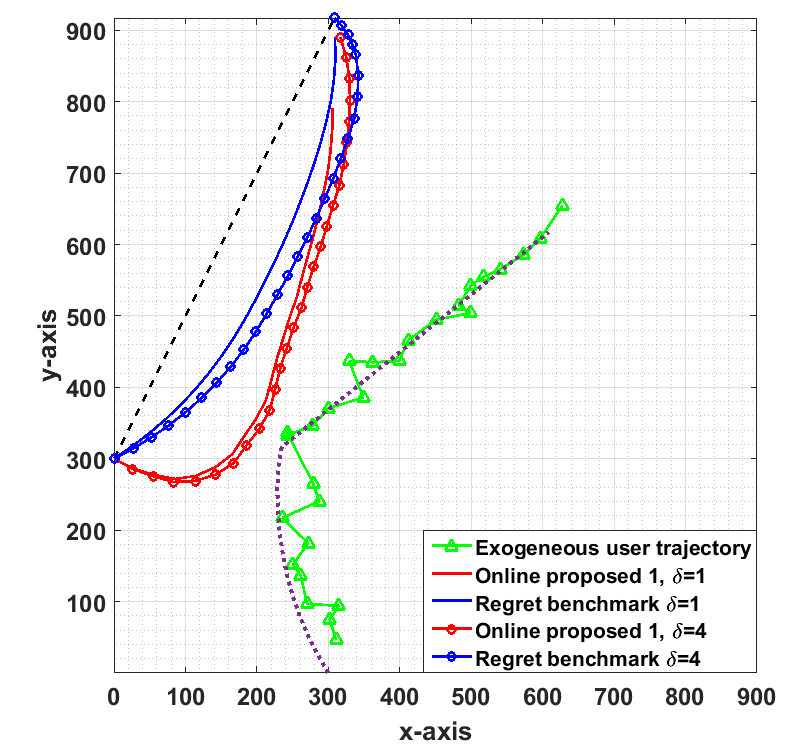} 
           \captionsetup{font=scriptsize}
           \caption{Proposed online method compared with regret benchmark of \eqref{online} for $\delta=1$ and $\delta=4$ (marked lines). \colb{Black dashed line is the straight line path shown for reference. The purple dotted line is the true trajectory of the exogenous user and the trajectory in green represents the  erroneous location estimates available to user $\x$.}}   
           \label{fig2}  
\end{figure}
\begin{figure}
 	\centering
    \includegraphics[scale=0.19]{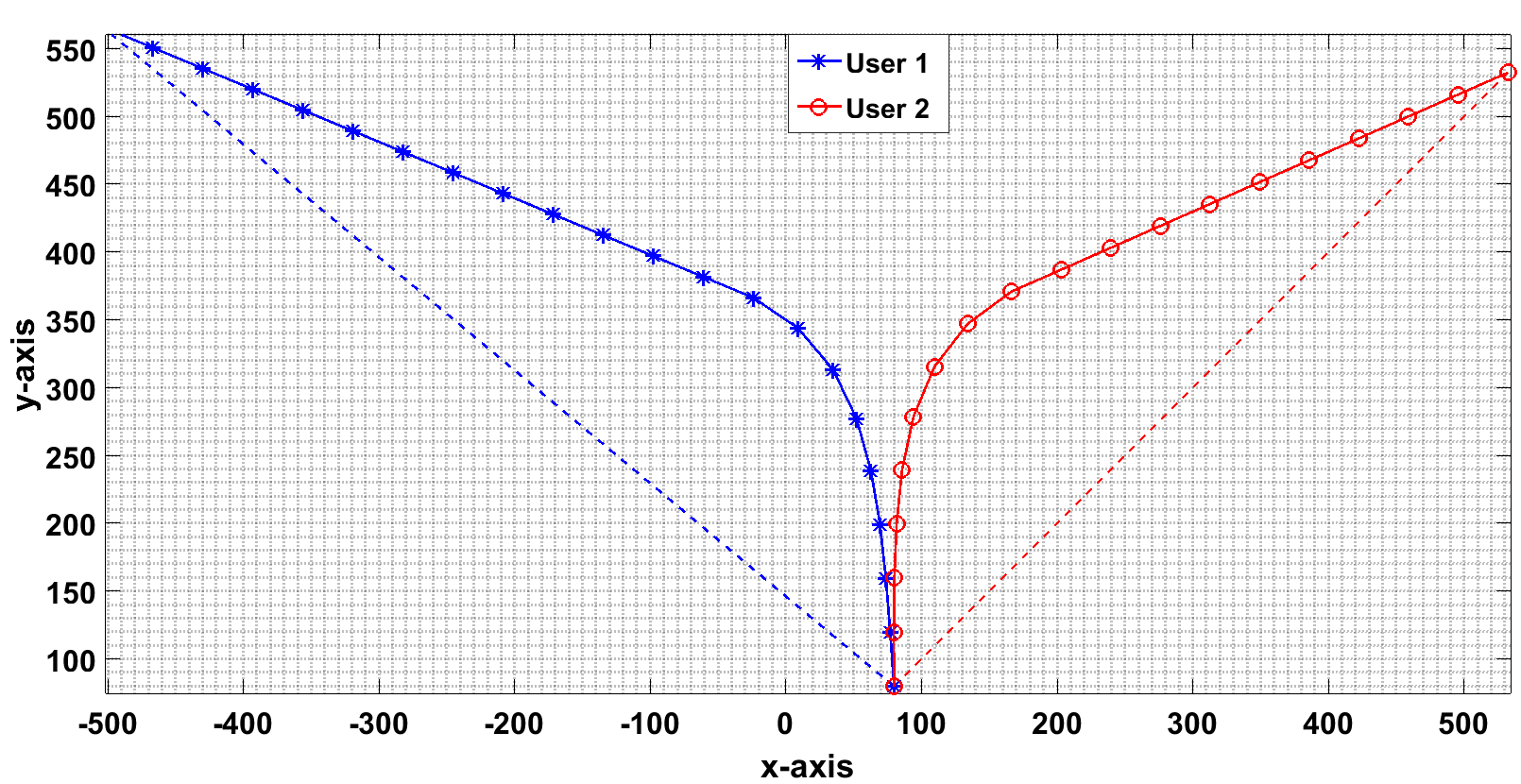}
 	\captionsetup{font=scriptsize}
 	\caption{D2D cooperative users with same starting and opposite destinations for $\delta=2$.  }
 	\label{second}
\end{figure} 
\begin{figure*}
	\centering
           \subfigure[]{\includegraphics[width=0.35\linewidth,height = 0.25\linewidth,trim={0cm 0cm 0cm 0cm},clip]{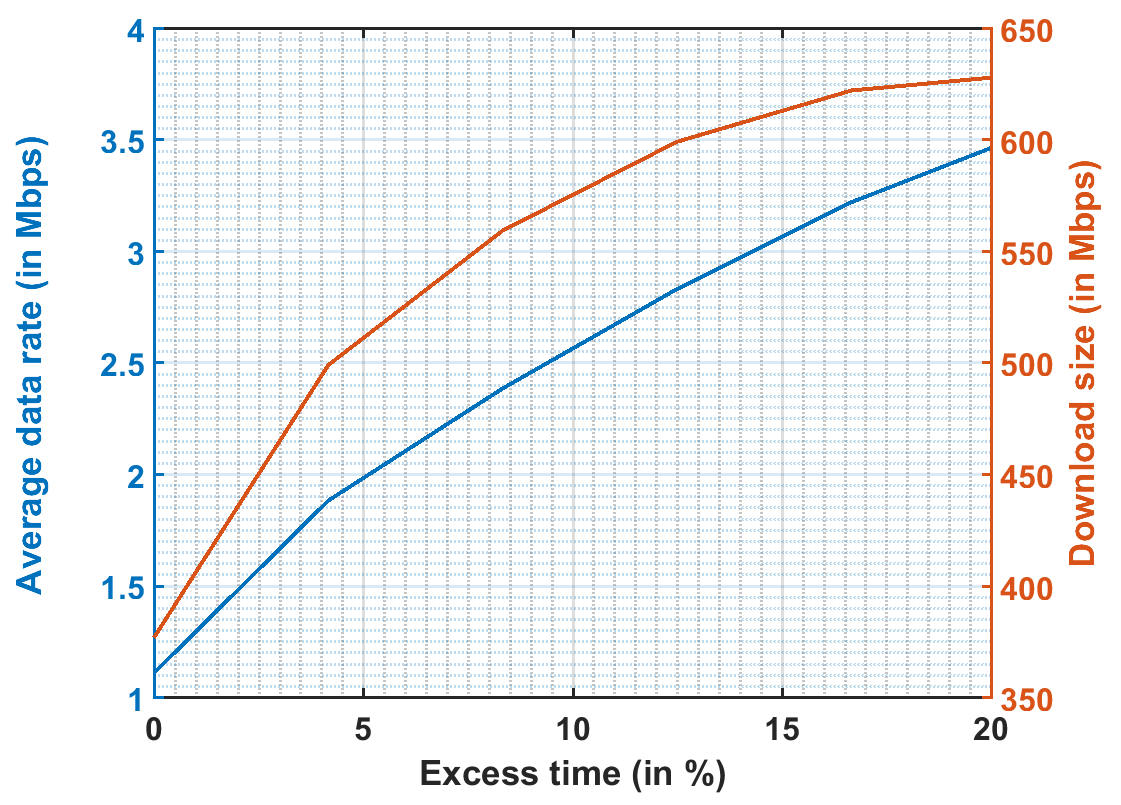}} 
           \hspace{-0.1cm}\subfigure[]{\includegraphics[width=0.3\linewidth,height = 0.25\linewidth,,trim={0cm 0cm 1cm 0cm},clip]{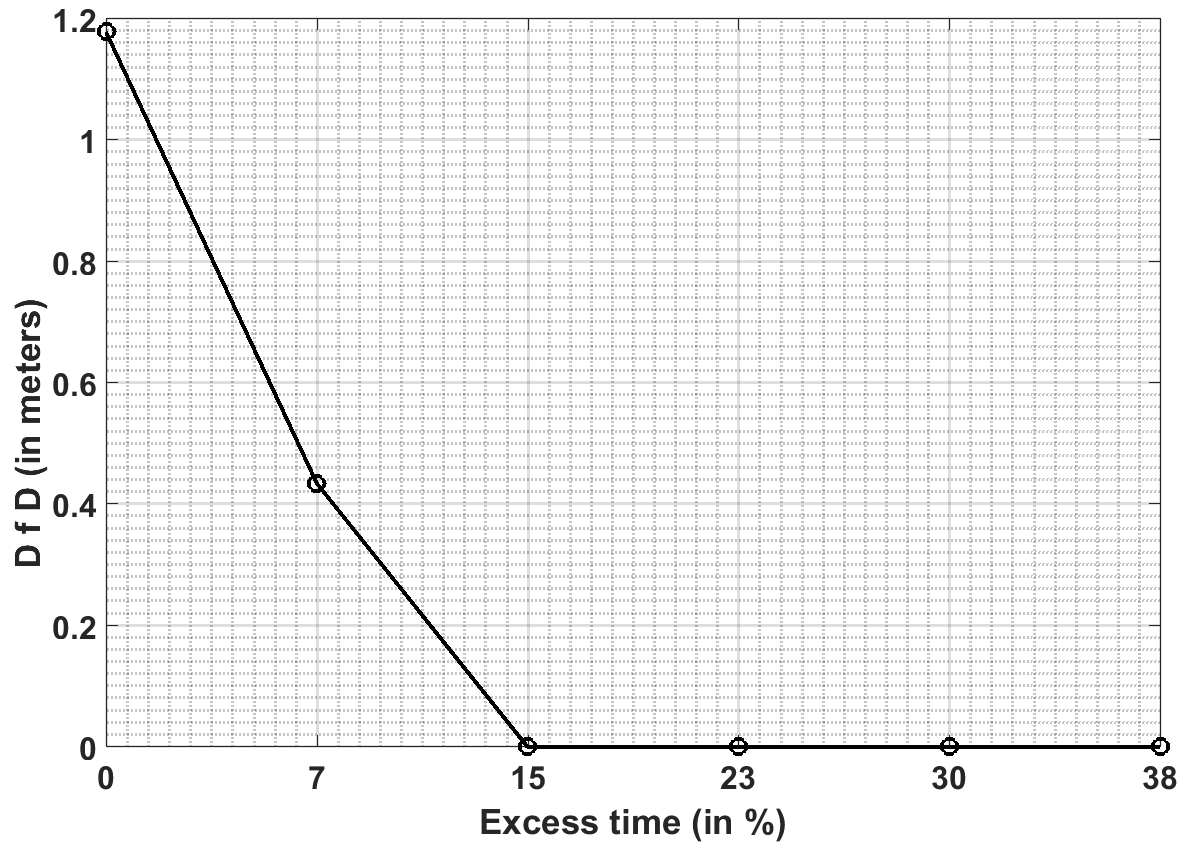}}  
           \subfigure[]{\includegraphics[width=0.3\linewidth,height = 0.25\linewidth,trim={0cm 0cm 0.1cm 0.1cm},clip]{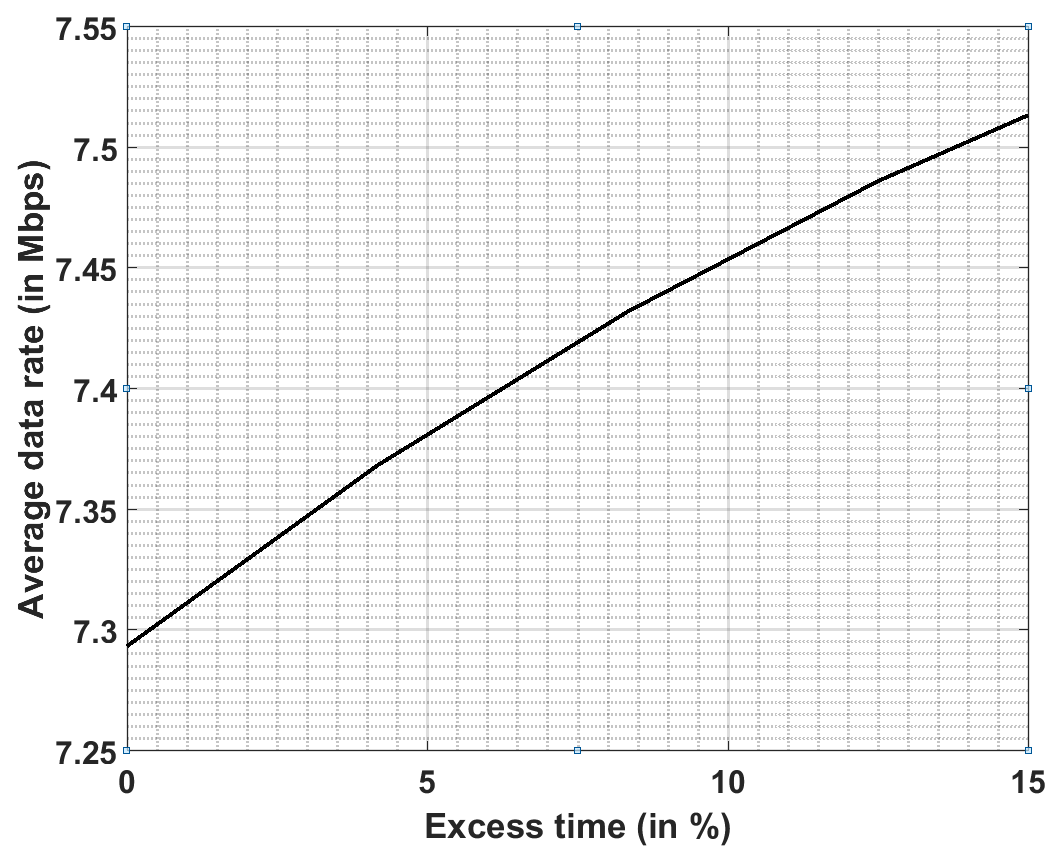}}
           \captionsetup{font=scriptsize}
           \caption{\textbf{(a)} Average data rate obtained with respect to increasing delay $\delta$ (excess time) for D2D cooperative users, and corresponding size of downloaded file. \textbf{(b)} Distance from the destination (D f D) vs percentage of excess time. \textbf{(c)} D2D non-cooperative user performance in terms of average data rate achieved versus excess time by using proposed online algorithm of \eqref{proposed}.} 
           \label{fig4}    
\end{figure*}
Finally, we solve for the optimal trajectory in an online fashion using online gradient ascent updates as shown in \eqref{online}. Specifically, for the objective function that takes the form as described in \eqref{huber} the OGA updates are given as follows
\begin{align}
&\hat{\x}{(t\!+\!1)} = P_\X(\hat{\x}(t) \!+\! \tfrac{\mu}{\gt}(\bell(t)\!-\!\hat{\x}(t)) \!+\! \tfrac{1\!-\!\mu}{\gt} P_v(\bell(t)-\hat{\x}(t)))\label{algo_main}.
\end{align}
A special case occurs when $\mu \approx 0$ and $\gamma = 1$, for which case, the updates become 
\begin{align}\label{proposed}
\hat{\x}(t+1) = \hat{\x}(t) + P_v(\bell(t)-\hat{\x}(t)).
\end{align}
Note that the gradient of the utility function is Lipschitz continuous with parameter $L = 1$ and strongly convex with parameter $\mu$. Moreover, if the user operates in a compact region of diameter $R$, the gradient is bounded as $\norm{\nabla U_t(\x)}_2 \leq \mu R + v(1-\mu)$. Consequently, (\textbf{A4}) is satisfied if $\mu R + v_{\max}(1-\mu) \leq \gt v_{\max}$ or equivalently, if $\gt \geq \frac{\mu (R-v_{\max})+v_{\max}}{v_{\max}} \geq 1+\frac{\mu R}{v_{\max}}$. For instance when $\mu$ is close to zero, we simply require $\gt \geq 1$. Next, we provide bounds to various problem parameters that evolve with the general constraint function $g_t(\x(t),\x(t+1)) := \norm{\x{(t+1)}-\x(t)} \leq \alpha(t) v_{\max}$ and recall that a special case is considered in all our experiments with $\alpha(t)= 1$.

\subsection{Choice of Step Size $\gamma(t)$ and Other Bounds}
\colb{In order to satisfy (A4), it is required that $\gt \geq \tfrac{\tilde{G}}{v(t)} = \tfrac{\tilde{G}}{v_{\max}\alpha(t)}$, where $\norm{\gradt U_t(\x(t))} \leq \tilde{G}$. Let us denote $\alfmin := \min_t\alpha(t)$ and $\alfmax = \max_t \alpha(t)$. Also recall that $\gmin := \min_t \gamma(t)$ and $\gmax := \max_t \gamma(t)$, so that $\gmax = \tfrac{\tilde{G}}{v_{\max}\alfmin}$ and $\gmin = \tfrac{\tilde{G}}{v_{\max}}$ where we set $\alfmax = 1$. Therefore the bounds in Theorem 1 of our manuscript (i.e., $L \leq \gmax \leq 2\gmin - L $) translate to
\begin{equation}
\tfrac{\tilde{G}}{2\tilde{G}-v_{\max} L } \ \  <  \ \ \alfmin \ \  <  \tfrac{\tilde{G}}{v_{\max}L}.
\end{equation}
Therefore, $\gt$ can be chosen as 
\begin{equation}
\text{choice of $\gt$} : \ \ \ \max\{\tfrac{\tilde{G}}{v_{\max}\alpha(t)},L\} < \gt < \tfrac{\tilde{G}}{v_{\max}\alfmin}.
\end{equation} 
In practice if the maximum norm gradient $\tilde{G}$ is not known a priori, one can use an estimate of the form $\max_{\tau\leq t} \norm{\gradt U_\tau(\x(\tau))}$.} Observe that the learning rate $\gamma(t)$ is indeed dependent on $t$, therefore achieving better utilities being adaptive to time-varying environmental conditions and constraints.

\subsection{Results for D2D Trajectory Planning Problem}  
This section provides detailed simulations for the various formulations provided here. Regardless of the algorithm used, the different formulations and settings will be compared from the cumulative D2D rates. At each time slot $t$, the maximum achievable rate given by 
\begin{align}
 R(t) = W\log \left(1+\frac{RSS}{RSS+\sigma^2}\right)\label{rate}.
\end{align}
Recall that $RSS=\norm{\x(t)-\y(t)}^{-\alpha_p}$ is the scaled path loss component and $\sigma^2$ is the appropriately scaled noise power. As remarked earlier the achievable rate does not represent the actual rate seen by the users, but is used here only for relative performance evaluation.

First, we present simulations for the online framework proposed in \eqref{proposed} for non-cooperative settings  where the users reveal their current locations but keep their future locations private. Here we consider $\mu=10^{-3}$ and $\gamma=1$. The online trajectory obtained for the proposed algorithm in  \eqref{proposed} is shown in Fig.~\ref{fig2} for different values of delay $\delta$. \colb{Note that in this case, in order to mimic a realistic scenario, we have considered a noisy trajectory of user~$\y$ while designing better data rate trajectory of user~$\x$. To this end, we have added Gaussian noise with zero mean and standard deviation of 1 meter}. In order to further examine the nature of the optimal trajectory, consider a special case when both users have the same starting point $\mathbf{s}_x = \mathbf{s}_y = [80\  80]^T$, but different destinations $\d_x = [-400\ 480]^T$ and $\d_y = [600 \ 600]^T $. The optimal trajectories are shown in Fig.~\ref{second} for $\delta = 2$. In this case, the trajectory exhibits a knee region, where the users initially stay close together and then abruptly split up to head towards their respective destinations. Note that, for all these cases, we assume the destination of user~$\x$ to be fixed but allow the trajectory of user~$\y$ to be exogenous.  As before we consider the square-law cost function and maximum velocity of $v_{\max} = 1$ m/s.

Next, consider the offline problem formulated in \eqref{main} for the case when $\X$ is simply a box in $\mathbb{R}^2$ and the utility function is $U(x) = -x^2$. We consider two pedestrians that start at the coordinates $\mathbf{s}_x = [0~400]^T$ and $\mathbf{s}_y = [400~0]^T$ and have destinations at $\d_x = [400~1200]^T$ and $\d_y = [800~800]^T$, respectively. Both users walk with the maximum speed of $v_{\max}^x = v_{\max}^y = 1$ m/s and take about $15$ minutes to reach their destinations through the direct path. In order to keep the problem sizes small, we divide this time into $T_x = T_y = 24$ time slots. As explained earlier, both the users are willing to incur an excess delay $\delta$ in reaching their destinations. To this end, the average D2D rate as a function of the allowable delay $\delta$ in this context is as shown in Fig.~\ref{fig4}(a). The second y-axis of this figure depicts the increase in the total downloaded file size with respect to delay.  The other parameters used in Fig.~\ref{fig4}(a) are $\alpha_p = 2.5$, $W = 10$ MHz, and $\sigma^2 = 0.2$. As expected, the average achievable rate continues to increase with $\delta$, as it allows larger deviations from the direct path. Next, we have compared the trajectories generated for various values of $\delta$ using the proposed IOGA against the offline regret benchmark defined in \eqref{online}. To plot this figure, the regularizer value used is $\lambda(t)=\frac{t}{T}$. If $\lambda(t)$ is not appropriately chosen, it may happen that the user will not be able to reach the destination. But this condition does not arise for sufficiently high $\delta$ as depicted in Fig.~\ref{fig4}(b). Note that for smaller values of $\delta$ (less excess time), the user is at a small distance from its final destination, which reduces to zero for larger values of delays (more excess time) and the user will always reach the destination. Further in Fig.~\ref{fig4}(c), the average data rate achieved by the proposed online algorithm is shown against the excess time.


\section{Energy-efficient trajectory planning in ocean environments}\label{energy_eff}
Consider a watercraft  operating in an ocean environment and seeking a possibly time-varying goal. While the watercraft has propulsive capabilities, it is required to reach the goal in a limited amount of time while expending minimal energy. The following information is available at any time $t$:
 \begin{itemize}
 	\item historical data on ocean currents in the region, e.g. available from Regional Ocean Modeling System (ROMS) \cite{roms},
 	\item current watercraft location $\xh(t)$,
 	\item current ocean velocity $\v_o(t)$ at location $\xh(t)$,
 	\item and current goal location $\d(t)$
 \end{itemize}
 in addition to problem parameters such as $T$, maximum watercraft velocity in still water $v_{\max}$, and other tunable parameters. 
 
 Designing online trajectories in such time-varying and uncertain environments is challenging and has never been attempted before. Existing methods usually rely on forecasts and require re-solving the full problem at every time slot. While such an approach can also be adopted here, we are interested in a more computationally efficient algorithm. Considering the uncertainty in ocean current measurements, the proposed algorithm will be based on the IOGA algorithm in \eqref{online} and therefore adheres to the guarantees in Theorem 1. However, in order to ensure assumption (A4) and to obtain reasonable performance, it is required to design the functions $U_t$ and $g_t$ as well as provide rules for choosing the step-size. 
 \begin{figure*}
 	     \centering
         \subfigure[]{\includegraphics[scale=0.25,trim={1.5cm 0.4cm 1.5cm 1.5cm},clip]{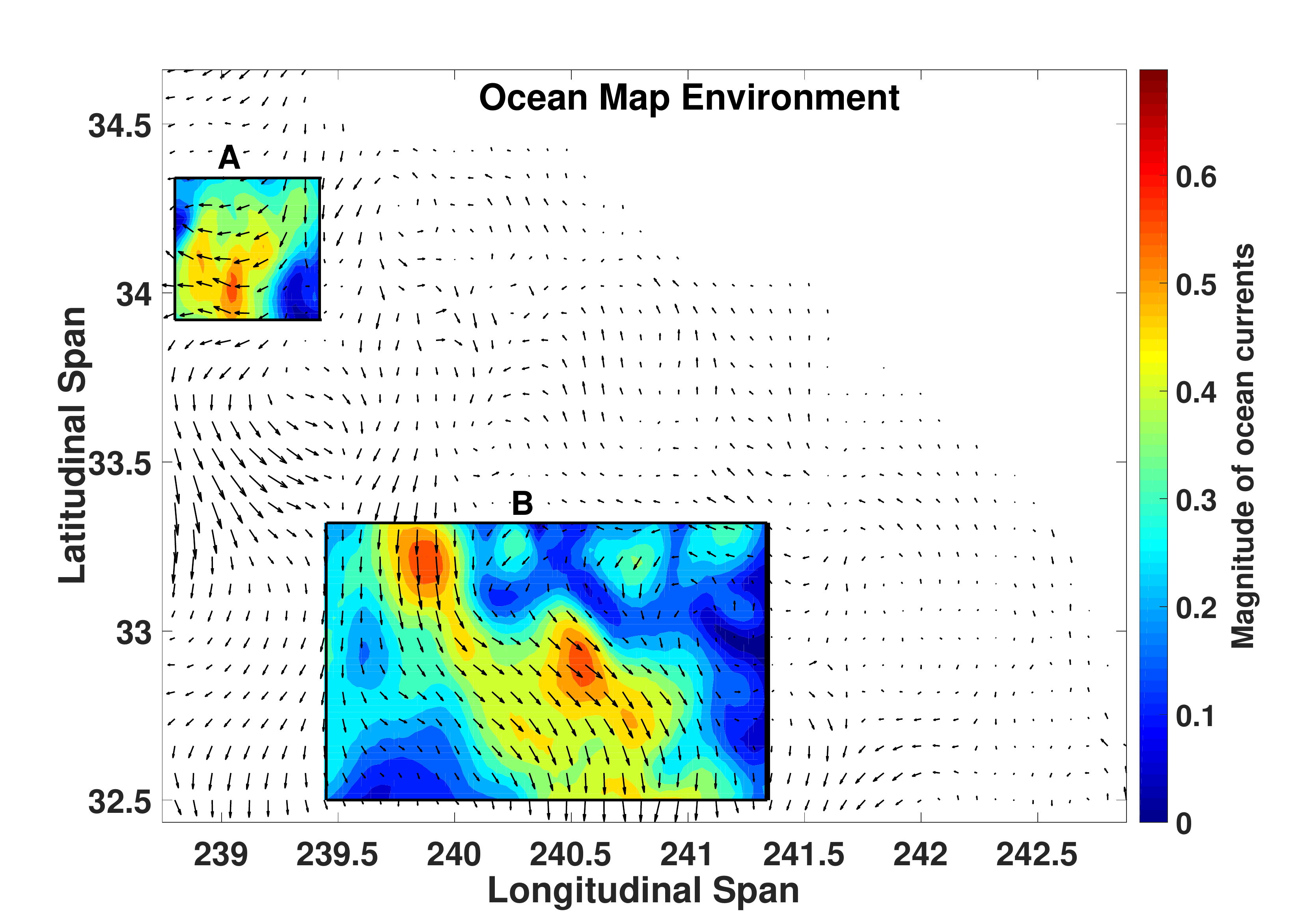}} \hspace{0.3cm}
         \subfigure[]{\includegraphics[scale=0.25, trim={1cm 0.27cm 1.5cm 0.6cm},clip]{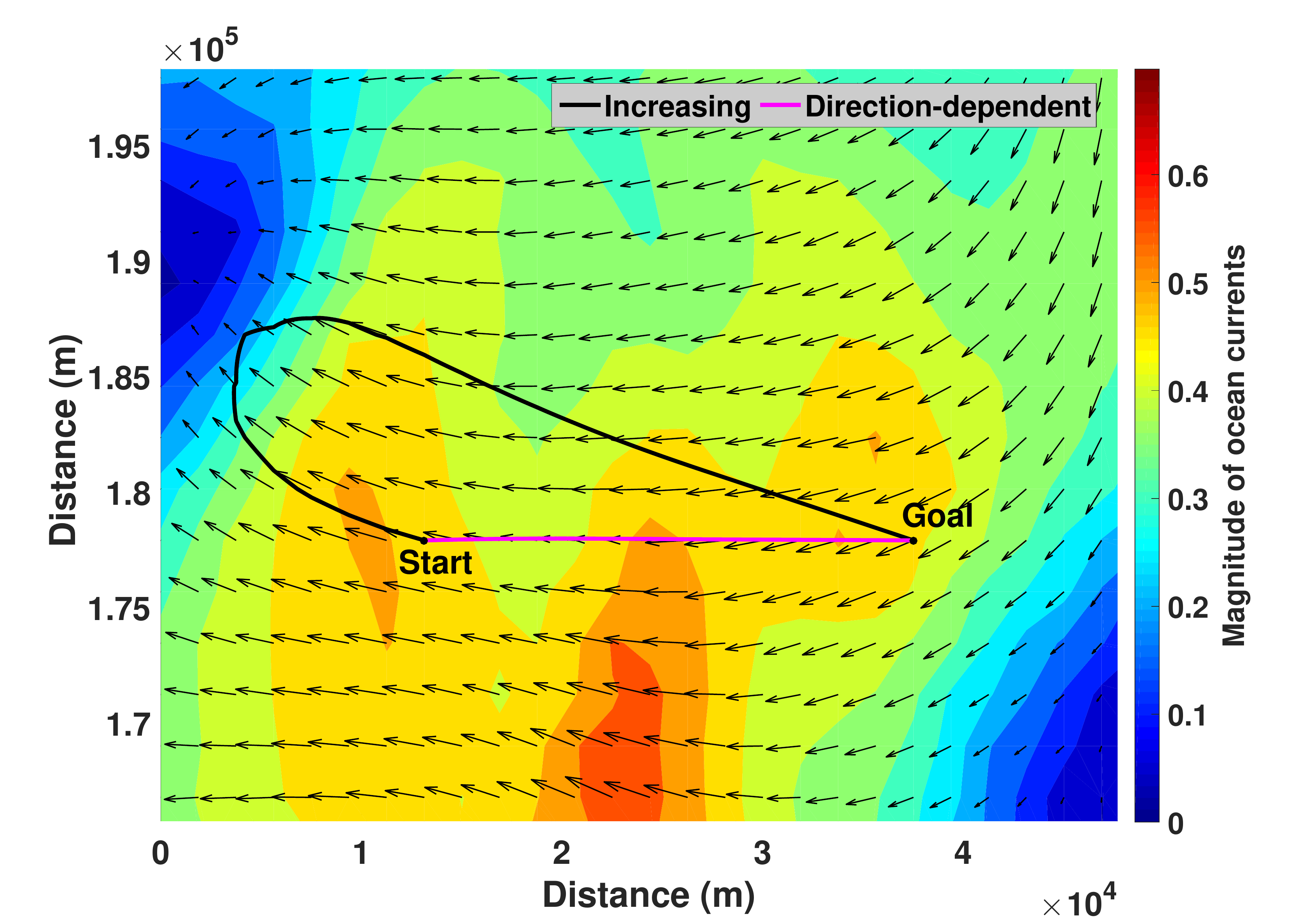}}\\
         \captionsetup{font=scriptsize}
         \caption{(a) Southern California Bight: Note that the ocean currents magnitudes at islands are made negligible (b) Demonstrating the behavior of Increasing vs. Direction-dependent strategies for $\lambda(t)$.}   
         \label{oceanmap}    
 \end{figure*}
 
\subsection{Design of the Utility Function $U_t$} We consider the following utility function 
\begin{align}\label{obj}
 	U_t(\x(t))  := & -\lambda(t) \norm{\x(t) - \d(t)}^2 \nonumber \\
 	&- (1-\lambda(t)) \ip{\xh(t-1)-\x(t), \v_o(t)}, 
 \end{align}
 where $\lambda(t) \in [0,1]$, is a control parameter. The utility function in \eqref{obj} is therefore a convex combination of the current squared distance from the target $\norm{\x - \d(t)}^2$ and the component of vehicle velocity in the negative direction of the ocean velocity. Consequently, an appropriate choice of $\lambda(t)$ ensures that the planned trajectory is headed towards the current goal location and in the direction of ocean currents, if possible. For instance, when $\lambda(t)$ is close to 1, the watercraft ignores the ocean currents and heads straight to the goal. Clearly, the choice of $\lambda(t)$ is critical towards ensuring the goal is reached in a timely manner and at the same time, energy consumption is minimized. Next, we discuss two different strategies for the choice of $\lambda(t)$. 
 
 \paragraph{Increasing $\lambda(t)$} The choice $\lambda(t) = t/T$ for $1\leq t\leq T$ is motivated from the observation that the watercraft may drift along the direction of the ocean currents initially to save energy. However, as time goes on, more importance must be placed on reaching the goal. While such a strategy was used in \cite{amrit}, it is agnostic to the direction and magnitude of $\v_o(t)$ and therefore suboptimal. In particular, the strategy fails when the ocean currents are pointing directly away from the goal, ultimately resulting in more energy expenditure than for the straight line path. Note that with this choice of $\lambda(t)$, Assumptions (A1)-(A3) are satisfied. 
 
 \paragraph{Direction-dependent $\lambda(t)$} Instead of increasing $\lambda(t)$, information about ocean directions may be taken into account using the choice:
 \begin{equation}\label{str2}
 \lambda(t) = 1 - \eta(t) \cos^2\left(\tfrac{\theta(t)}{2}\right),
 \end{equation}
 where $\theta(t) = \angle\left(\d(t)-\xh(t),\v_o\right)$, $\eta(t) = \frac{\norm{{\v_o(t)}}}{v_o^{\max}}$, and $v_o^{\max}$ is the maximum ocean current velocity estimated from historical data. Here, the ratio $\frac{\norm{{\v_o(t)}}}{v_o^{\max}}$ is the relative strength of the ocean currents at a given point. When the ocean currents are weak, i.e., $\norm{{\v_o(t)}} \ll v_o^{\max}$ or when the currents are directed away from the goal location, i.e., $\theta \approx \pi$ radians, $\lambda(t)$ is close to 1 and consequently the watercraft is headed towards the target. However, when the ocean currents are in a favorable direction, we have that $\lambda(t) < 1$, and therefore the watercraft may utilize the ocean velocity in order to save energy. Table \ref{table1} provides some example values of $\lambda(t)$ for various values of $(\eta(t),\theta(t))$. 

 \begin{table}[] \fontsize{8pt}{12pt}\selectfont
 \centering
  \captionsetup{font=scriptsize}
 	\caption{Values of $\lambda(t)$ for different $(\eta(t),\theta(t))$}
 	\label{table1}
 	\begin{tabular}{|c|c|c|c|}
 		\hline
 		\multirow{2}{*}{\begin{tabular}[c]{@{}c@{}}Direction of \\ ocean currents \\ w.r.t goal \end{tabular}} & \multicolumn{3}{c|}{Magnitude of ocean currents} \\ \cline{2-4} 
 		& \begin{tabular}[c]{@{}c@{}}Strong \\ $\eta(t) \approx 1 $\end{tabular} & \begin{tabular}[c]{@{}c@{}}Moderate \\ $\eta(t) \approx 0.5 $\end{tabular} & \begin{tabular}[c]{@{}c@{}}Lower \\ $ \eta(t) \approx 0 $\end{tabular} \\ \hline
 		$\theta(t) \in [0,\frac{\pi}{2}) $ & $\lambda(t) \in [0 , 0.5)$ & $\lambda(t) \in [0.5 , 0.75)$ & $\lambda(t) \approx 1$ \\ \hline
 		$\theta(t) \in [\frac{\pi}{2},\pi] $ & $\lambda(t) \in [0.5 , 1]$ & $\lambda(t) \in [0.75 , 1)$ & $\lambda(t) \approx 1$ \\ \hline
 	\end{tabular}
 \end{table}

 \subsection{Design of the Constraint Function $g_t$} Recalling that the average velocity of the watercraft at time slot $t$ is $\x(t+1)-\x(t)$ and measured in meters per slot, we consider the constraint function of the form
 \begin{align}\label{constr}
 	g_t(\x(t+1),\x(t)) := &\norm{\x(t+1)-\x(t)-\v_o(t)} - \alpha(t) v_{\max}, 
 \end{align}
 where $\left(\x(t+1)-\x(t)-\v_o(t) \right)$ represents the velocity of the watercraft relative to the ocean. Note that the magnitude of the relative velocity of the watercraft is physically limited to at most $v_{\max}$ meters per slot. The additional factor of $\alpha(t) \in (0,1]$ is included to restrict the maximum relative velocity further if required. When the objective function with $\lambda(t)$ as in \eqref{str2} is used, it only provides directional information and $\alpha(t)$ must be carefully chosen to ensure that assumption (A4) is satisfied. To this end, we consider
 \begin{equation}\label{alpha}
 \alpha(t) = \exp\left(-\beta \left(\frac{\delta}{T} + \eta(t) \cos\left(\frac{\theta(t)}{2}\right) \right) \right),
 \end{equation}
 where $\delta$, $\eta(t)$, and $\theta(t)$ are as defined earlier. Intuitively, $\alpha(t)$ is small when the ocean currents are in a favorable direction ($\theta \ll \pi/2$) and are either strong (large $\eta(t)$) or sufficient excess time is available (large $\delta$) to reach the goal. In such a scenario, the watercraft simply drift along the currents and can conserve energy. On the other hand, the watercraft runs the engine at full capacity when $\alpha(t)$ is close to one, such as when the ocean currents are in the opposite direction, and no delay is allowed ($\delta \approx 0$). The tuning parameter $\beta$ is dependent on the environment and must be learned a priori using the forecast data. Finally, the IOGA updates are applied using the objective function in \eqref{obj} while the step-size $\gt$ is chosen to satisfy Assumption (A4) where $g_t$ is chosen as per \eqref{constr}. 
\subsection{Choice of Step Size $\gt$ and Other Bounds}
\colb{Now recalling the constraint function defined under the energy efficient trajectory planning application, the IOGA updates results in quadratic inequality in $\gt$ 
\begin{equation}
\norm{\frac{1}{\gt} \gradt U_t(\x(t)) - \v_o(t) }_2^2 \leq \alpha(t)^2\ v_{\max}^2. 
\end{equation} 
Following the IOGA updates results in quadratic inequality in $\gt$. To this end, we solve it for $\gt$ with equality, pertaining to the application as follows  
\begin{align}\label{eqre}
\left(\norm{\v_o}_2^2  - \alpha(t)^2 v_{\max}^2 \right) \gt^2 &- \left(2 \gradt U_t(\x(t))^T \v_o(t)\right) \gt \nonumber \\ 
& + \norm{\gradt U_t(\x(t))}_2^2 \  = \ 0.
\end{align}}  	

\colb{For atleast one of the roots of equation \eqref{eqre} to be positive it is required that 
\begin{equation}
\alpha(t)\ \  > \ \  \frac{\norm{\v_o(t)}}{v_{\max}}, 
\end{equation}
holding for all $t$. Considering this lower bound one can define $\gmin$ and $\gmax$ from \eqref{eqre} as a function of $\alfmin$. 
Now using the lower bound on $\gmax > L $, we get     
\begin{equation}
\alfmin \ \ < \ \ \sqrt{\frac{\beta' + L\norm{\v_o(t)}^2}{L v_{\max}^2}}  \ \ < \ \ \sqrt{\frac{\beta' + L (v_o^{\max})^2}{L v_{\max}^2}}, 
\end{equation}
where $\alfmin = \min_t \left\{\alpha(t) \right\}_{t=1}^T$ and 
\begin{align}
\beta' \ = \ & \sqrt{\left(\begin{smallmatrix}
	v_{\max}^2\norm{\gradt U_t(\x(t))}^2 -2( \norm{\v_o(t)}^2 \norm{\gradt U_t(\x(t))}^2 \\ 
	-(\gradt U_t(\x(t))^T \v_o(t))^2 )
	\end{smallmatrix}\right)} \nonumber \\
& + \ \gradt U_t(\x(t))^T\v_o(t). 
\end{align}
Note that the lower bound on $\alfmin$ can be found by using the inequality $\gmax < 2\gmin -L$, that requires solving a cubic equation in terms of $\alfmin$. }

\subsection{Results for Energy-Efficient Trajectory Planning Problem} 
We tested our algorithm in a real oceanic environment built from historical ROMS ocean currents data taken from \cite{roms}. Specifically, we have chosen Southern California Bight output collected on August $10$, 2018, as  shown in Fig. \ref{oceanmap}(a).  
From the dataset, it is found that the ocean currents' velocity is found to be varying from $0.001\ m/s$ to a maximum of $0.69\ m/s$. Thereby we assume that the maximum velocity of the watercraft in still water is $1 m/s$ to make sure that it can travel forward even in strong disturbances. The energy incurred for travelling in planned trajectories is calculated in the same manner as detailed in \cite{witt} and is given by 
\begin{equation}
E = c_d V_r^3t
\end{equation}
where $c_d$ is the vehicles drag coefficient, $V_r$ is the magnitude of the required velocity for the motors to provide and $t$ is the travel time for the relevant section of the path.

\subsubsection{{Static goal}}
\colb{We first illustrate the working of the algorithm on a simple setting in order to understand the effect of the different choices of $\lambda(t)$. To this end, we consider a static goal scenario under noiseless gradient feedback.} Consider a part of the dataset (labeled `A' in Fig. \ref{oceanmap}(a)) where the ocean currents with high magnitude are directed away from the goal. The start and goal locations are as specified in Fig. \ref{oceanmap}(b). It can be observed that the IOGA algorithm with increasing $\lambda(t)$ makes the USV go away from the goal initially when $\lambda(t)$ is small. As time goes on, $\lambda(t)$ increases and more emphasis is placed on reaching the goal, and the USV again starts going in the direction of the goal. Such a trajectory is suboptimal since it requires the USV to travel a long distance against the current. In contrast, the direction-dependent choice of $\lambda(t)$ makes the USV travel directly to the goal and incurs less energy, as also evident from Table \ref{table2}. 
\begin{figure} 
      \centering  
      \subfigure{\includegraphics[width=0.98\linewidth,trim={0.2cm 0cm 0cm 0.1cm},clip]{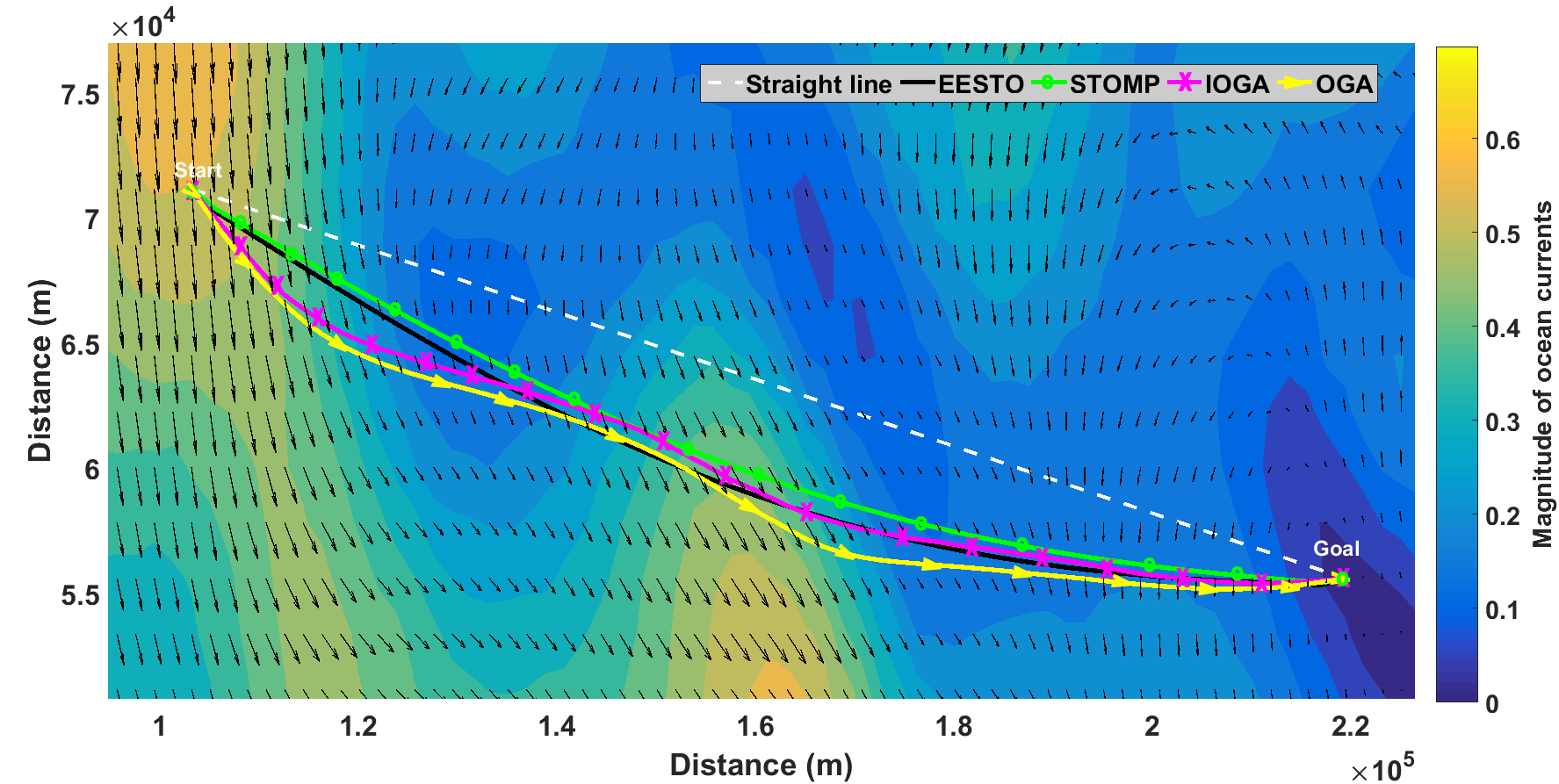}} 
      \captionsetup{font=scriptsize}
      \caption{Trajectory comparisons of STMOP \cite{stomp}, EESTO\cite{eesto} and IOGA with energy efficient control strategy. Note that straight line trajectory with broken lines is shown for reference.l }\label{traj}
\end{figure}  
\begin{table}[ht]
\captionsetup{font=scriptsize}
\caption{Energy cost comparisons}
\label{table2}
\centering
\begin{tabular}{@{}m{2.5cm}m{2.5cm}@{}}
\toprule[1 pt] 
\textbf{$\lambda(t)$} &  Energy cost (kJ) \\
\hline
\vspace{2mm}
\textit{Increasing }& \hspace{3mm}184.42\\
\vspace{2mm}
\textit{Direction-dependent} & \hspace{3mm}\textbf{130.29}  \\
\bottomrule[1pt]  
\end{tabular}
\end{table}

Next, to illustrate the performance of the proposed IOGA algorithm in comparison to the state-of-the-art algorithms under a static goal scenario, we now choose another chunk of real data `B' as shown in Fig. \ref{oceanmap}(a). \colb{The start and goal locations are as specified in the Fig. \ref{traj} and are separated by around 459 km. For the chunk of the data chosen, we perturb the ocean current forecast by adding zero mean white Gaussian noise with standard deviation varying from $0\%$ to $10\%$ of the maximum ocean current disturbance recorded. This is done to consider the uncertainty in the ocean currents while assuming the available forecast data as the ground truth. Trajectories are generated using STOMP \cite{stomp}, EESTO \cite{eesto} and the proposed IOGA with the direction-dependent control strategy.} 
	
\colb{In Fig. \ref{traj} we chose to display trajectories for a particular case where the standard deviation of the perturbation in ocean currents is $5\%$. The figure also includes the trajectory designed using IOGA with true gradient feedback (i.e., OGA algorithm) for reference along with the true ocean current forecast.}  Since STOMP \cite{stomp} and EESTO \cite{eesto} are sampling based methods they do not output the same trajectory at every run. Henceforth, we executed those algorithms 100 times each and picked the trajectory with the minimum energy cost incurred. 

\begin{figure} 
      \centering  
      \subfigure{\includegraphics[width=0.9\linewidth,trim={0cm 0cm 0cm 0cm},clip]{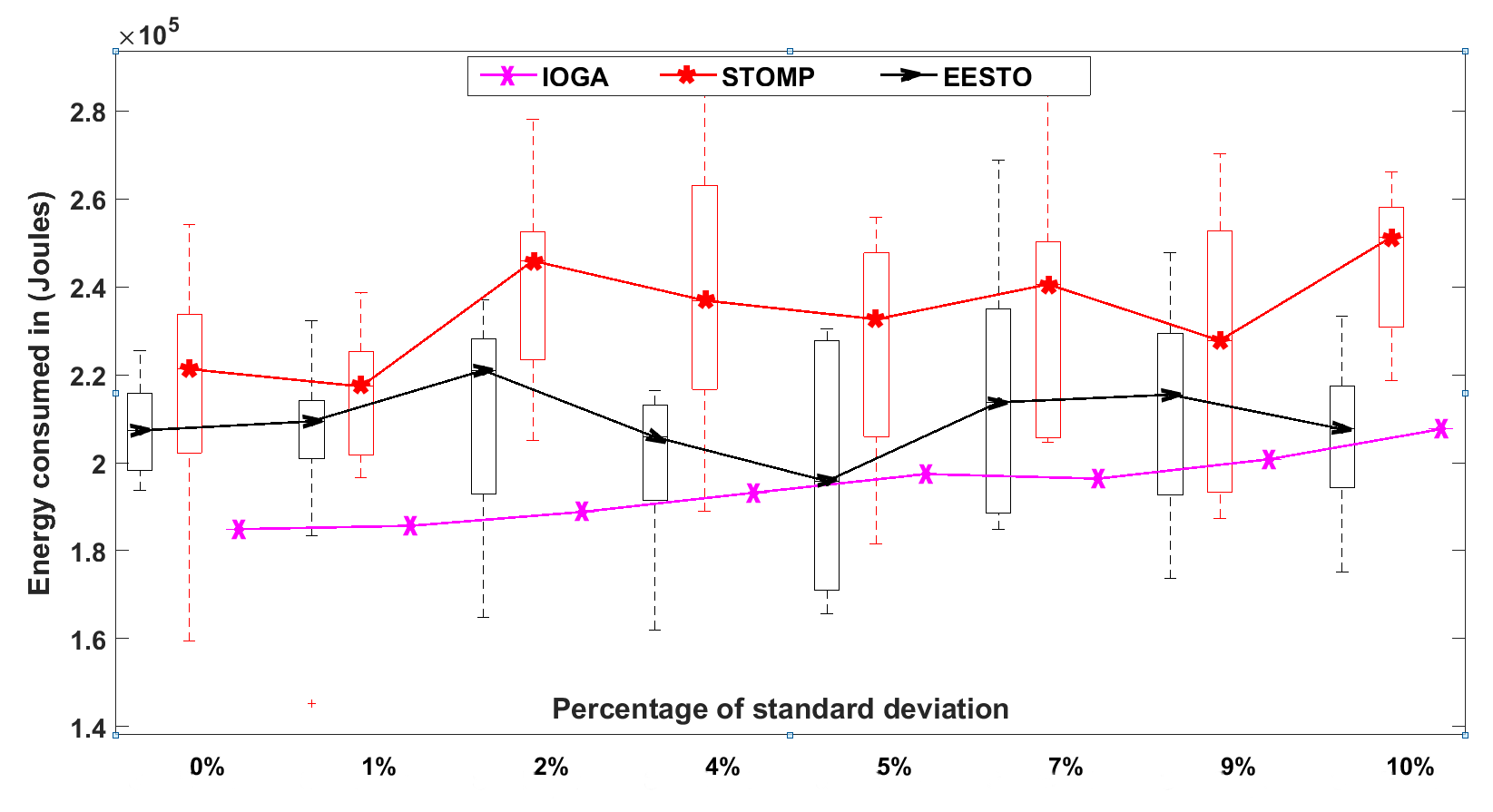}}
      	\captionsetup{font=scriptsize}
      \caption{\colb{Energy cost comparisons for STOMP\cite{stomp}, EESTO \cite{eesto}, IOGA against various noise levels}}\label{boxp}
\end{figure}

The variation of the total energy consumed by all the three algorithms against varied perturbations of ocean currents is shown in Fig. \ref{boxp}. It is evident from Fig. \ref{boxp} that the proposed algorithm is generating trajectories that are more energy efficient than the state-of-the-art methods besides being online. More importantly, unlike STOMP\cite{stomp} and EESTO \cite{eesto}, the obtained trajectory is the same for each run, and the resulting output has zero variability. 

\begin{figure} 
      \centering  
      \subfigure{\includegraphics[width=0.98\linewidth,trim={0cm 0cm 0cm 1.2cm},clip]{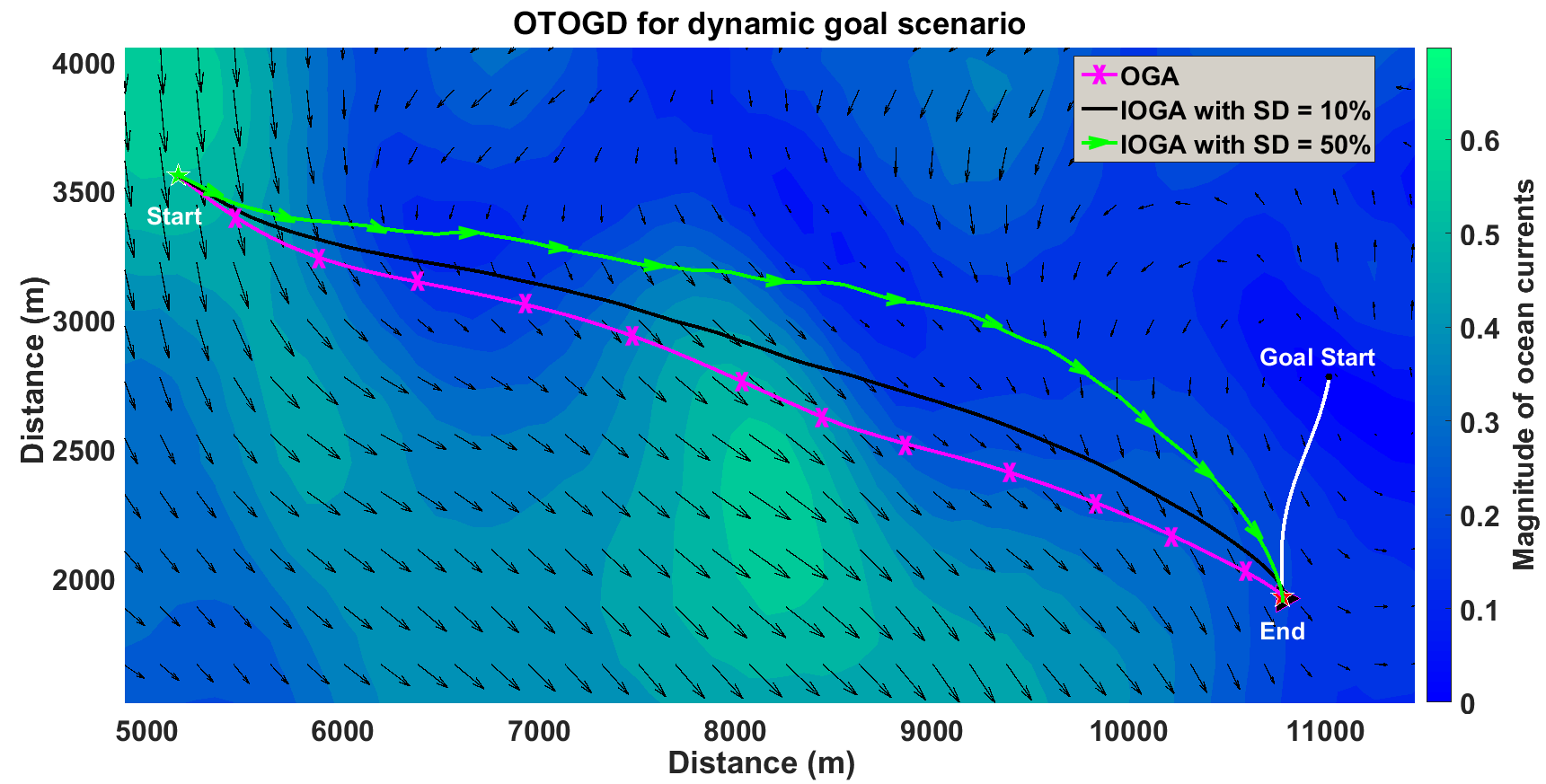}}
      	\captionsetup{font=scriptsize}
      \caption{Illustration of planning in moving goal scenario. {We display trajectories generated by IOGA algorithm for various levels of noise characteristics i.e., standard deviation (SD) = \{0\%, 10\%, 50\%\} of the maximum ocean current velocity magnitude.}   }\label{moving_goal}
\end{figure} 
 
\subsubsection{{Dynamic goal}}
To illustrate the performance of the proposed algorithm under dynamic goal scenario, we make use of the same chunk `B' of the real data. Here we simulated a moving goal scenario and assumed that the exact location of the goal at every time instant is known to the agent. Note that the proposed framework still requires the maximum velocity of the goal to be significantly less than that of the USV. \colb{Fig. \ref{moving_goal} displays the final trajectories for various cases of IOGA algorithm along with the trajectory traversed by the goal. Specifically we display trajectories generated by IOGA algorithm for various levels of uncertainties in the ocean current measurements. Observe the variation in the shape of the trajectories displayed. In the regions where ocean current magnitude is significantly low, the algorithm ignores the effect of ocean currents and generates trajectories in the direction of goal.}

It is remarked that in theory, existing algorithms such as \cite{eesto} can perform re-planning at every time instant thereby following a moving goal. However, with per-waypoint re-planning, EESTO incurs a computational cost of $\mathcal{O}(IKN^4)$ where $K$ is the number of samples, $N$ represents the total number of waypoints, and $I$ is the total number of iterations required for convergence. In contrast, the proposed algorithm only incurs a cost of $\mathcal{O}(N)$ even when the goal and the environment are time-varying.  
\subsection{Delay vs Energy Conservation Analysis}
To study the importance of trade off between delay and the amount of energy conserved, we ran STOMP \cite{stomp}, EESTO \cite{eesto} along with the proposed IOGA algorithm with noiseless gradient feedback for 10\% and 30\% allowed delay in time for 30 randomly chosen chunks of real data.
For each chunk, Table \ref{table3} shows the average amount of energy conserved as compared to the energy incurred in traveling along the straight line path from the starting point to the goal. \colb{ Note that, we calculate the amount of energy conserved by taking the total energy incurred in traversing the straight line path as reference.} As evident from the table \ref{table3}, the proposed algorithm saves more energy as compared to state-of-the-art algorithms. Further, existing algorithms exhibit a high variance in the conserved energy value, requiring multiple runs.  
 
\begin{table}[]\fontsize{8pt}{12pt}\selectfont
	\captionsetup{font=scriptsize}
\caption{Energy conservation vs. delay analysis}
\label{table3}
\begin{tabular}{|c|c|c|c|}
\hline
\multirow{2}{*}{\begin{tabular}[c]{@{}c@{}}$\%$ of \\ excess time slots\end{tabular}} & \multicolumn{3}{l|}{\hspace{1cm} Energy conserved (kJ)} \\ \cline{2-4} 
 & EESTO & STOMP & IOGA \\ \hline
10$\%$ & 112.01 $\pm$ 24.4 & 107.10 $\pm$ 28.2 & \textbf{123.08} \\ \hline
30$\% $ & 184 $\pm$26.6 & 170.2 $\pm$ 31.1 & \textbf{217.13} \\ \hline
\end{tabular} 
\end{table}


\section{Conclusion} 
This paper considered the utility optimal trajectory design problem with time-varying coupling constraint functions. The problem is formulated within convex optimization formalism and solved by employing a modified gradient ascent algorithm under noisy gradient feedback. We have referred to this algorithm as inexact online gradient ascent algorithm (IOGA). Next, we have shown that the proposed IOGA algorithm incurs sublinear offline regret under certain mild assumptions. Moreover, assuming that the reliable predictions of the environment were available, we proposed a modified IOGA update rule which is shown to incur improved sublinear offline regret as opposed to the standard update rule. More importantly, a lower bound to the offline dynamic regret is derived that helps in judging the optimality of any online policy developed in this context. 

The generality of the proposed framework helped in applying the algorithm to two varied applications (i) planning data rate maximizing trajectories for commute users in D2D communication setting (ii) planning energy efficient trajectories under strong and uncertain oceanic disturbances \footnote{MATLAB implementation is made available at https://tinyurl.com/traj-planner}. The proposed algorithm is shown to handle real world uncertainties without failing catastrophically and yields utility optimal trajectories as compared to the state-of-the-art algorithms without relying much on the forecast data. However, the present formulation cannot directly incorporate roads in (D2D application) and obstacles in (ocean environments), since the bounds developed here rely on the area being convex. Another open issue that is relatively common is the case of multi-agent planning, possibly interacting
during different parts of the commute. For instance in the D2D application user 1
may communicate with user 2 for the first 10 minutes and subsequently switch to a new user 3.

\ifCLASSOPTIONcaptionsoff
\fi
\bibliographystyle{IEEEtran}
\begin{footnotesize}
\bibliography{refer}
\end{footnotesize}

 \newpage
\onecolumn 
\section*{Supplementary Material}
\subsection{Proof of Theorem \ref{thm}} \label{app1}


\begin{IEEEproof}
	Since $-U_{t+1}$ is $\mu$-convex from Assumption {(A1)}, we have that 
	\begin{align}\label{sc}
	U_{t+1}(\x^\star(t+1)) \leq & U_{t+1}(\xh(t)) - \tfrac{\mu}{2} \norm{\x^\star(t+1) - \xh(t)}_2^2 \ip{\nabla U_{t+1}(\xh(t)),\x^\star(t+1) - \xh(t)},
	\end{align}
	for  all $1\leq t \leq T-1$. Adding and subtracting $\ip{\nabla U_{t+1}(\xh(t)), \xh(t+1)} $ to the right of \eqref{sc} and rearranging, we obtain 
	\begin{align}
	U_{t+1}(\x^\star(t+1)) + \tfrac{\mu}{2} \norm{\x^\star(t+1) - \xh(t)}_2^2  &\leq \  U_{t+1}(\xh(t)) + \ip{\nabla U_{t+1}(\xh(t)),\x^\star(t+1)-\xh(t+1)} \nonumber \\  &\ \ \ + \ip{\nabla U_{t+1}(\xh(t)),\xh(t+1) - \xh(t)}, \label{utdiff}
	\end{align}
	for $1\leq t \leq T-1$. Next observe that the IOGA update in \eqref{online} can be written as 
	\begin{align}
	\!\!\!\!\xh(t\!\!+\!\!1) \!= \!\arg\max_{\x\in\X} &\ip{{\nabla\Utt(\xh(t))}, \x\!-\!\xh(t)} \!-\! \tfrac{\gt}{2}\norm{\x\!-\!\xh(t)}^2_2, \label{update}
	\end{align}
	whose optimality condition implies that
	\begin{align}
	\!\!\!\!\!\!\ip{\nabla {\Utt(\xh(t))} \!-\! \gt  (\xh(t\!\!+\!\!1)-\xh(t)),\x^\star(t\!\!+\!\!1)\!-\!\xh(t\!\!+\!\!1)} \!\leq\! 0, \label{optcon1}
	\end{align}
	since $\x^\star(t+1)\in\X$.  Recall that { $\nabla \Utt(\xh(t)) = \nabla {U}_t(\xh(t)) + \nt$} for some $\nt \in \R^2$. Let us introduce $G_t(\xh(t)) :=  \nabla U_{t+1}(\xh(t)) - \nabla U_t(\xh(t))$ so that the optimality condition in \eqref{optcon1} can be written as
	\begin{align}
	\ip{\nabla  U_{t+1}(\xh(t)),\x^\star(t+1)-\xh(t+1)} &\leq \gt \ip{\xh(t+1)-\xh(t),\x^\star(t+1)-\xh(t)} + \ip{G_t(\xh(t)),\x^\star(t+1) - \xh(t+1)}  \nonumber \\ 
	&\ \ \ -\gt\norm{\xh(t+1) - \xh(t)}_2^2 {- \ip{\nt,\x^\star(t+1) - \xh(t+1)}},  \label{int2} 
	\end{align}
	which upon substituting into \eqref{utdiff} yields for $1\leq t \leq T-1$
	\begin{align}\label{scf}
	U_{t+1}(\x^\star(t+1)) + \frac{\mu}{2} \norm{\x^\star(t+1) - \xh(t)}_2^2  
	& \leq U_{t+1}(\xh(t)) + \ip{\nabla U_{t+1}(\xh(t)),\xh(t+1) - \xh(t)} + \ip{G_t(\xh(t)),\x^\star(t+1) - \xh(t+1)}  \nonumber \\ 
	&\ \ \ - \gt \norm{\xh(t+1) - \xh(t)}_2^2 {- \ip{\nt,\x^\star(t+1) - \xh(t+1)}} \nonumber \\
	&\ \ \ + \gt\ip{\xh(t+1)- \xh(t), \x^\star(t+1) - \xh(t)}. 
	\end{align}
	Since $U_{t+1}$ is $L$-smooth from Assumption {(A2)}, the quadratic lower bound implies that 
	\begin{align}\label{lp}
	U_{t+1}(\xh(t+1)) \geq U_{t+1}(\xh(t)) - \tfrac{L}{2}\norm{\xh(t+1)-\xh(t)}_2^2  + \ip{\nabla U_{t+1}(\xh(t)),\xh(t+1)-\xh(t)},
	\end{align}
	for $1\leq t \leq T-1$. Combining \eqref{scf} and \eqref{lp} we obtain 
	\begin{align}\label{eq11}
	U_{t+1}(\x^\star(t+1)) - U_{t+1}(\xh(t+1)) & \leq \left(\tfrac{L}{2} - \gamma(t) \right)\norm{\xh(t+1) - \xh(t)}_2^2 - \tfrac{\mu}{2}\norm{\x^\star(t+1) - \xh(t)}_2^2  \nonumber \\
	&\ \ \ + \gt \ip{\xh(t+1) - \xh(t), \x^\star(t+1) - \xh(t)}  + \ip{G_t(\xh(t)),\x^\star(t+1) - \xh(t+1)} \nonumber \\
	&\ \ \ \ {- \ip{\nt,\x^\star(t+1) - \xh(t+1)}},
	\end{align}
	for $1\leq t \leq T-1$. Note that since  $\{\x^\star_t\}_{t=1}^T$ is a solution to the off-line problem, it satisfies
	\begin{align} \label{wrg}
	\s U_{t}(\x^\star(t)) \geq \s U_{t}(\y(t)),
	\end{align}
	for any $\{\y(t)\in\X\}_{t=1}^T$ satisfying the time varying coupling constraints. For instance, the sequence $\{\xh(t) \in \X\}_{t=1}^T$ satisfies $g_t(\x(t),\x(t+1)) \leq 0$ for all $1\leq t \leq T-1$ from assumption ({A4}). Note that the first term on both sides of \eqref{wrg} is the same since $\x^\star(1) = \xh(1) = \mathbf{s}$ and can be dropped. Equivalently, the condition in \eqref{wrg} for $\y(t) = \xh(t)$ can be written as
	\begin{equation}\label{optim}
	\sum_{t=1}^{T-1} U_{t+1}(\x^\star(t+1)) \geq \sum_{t=1}^{T-1} U_{t+1}(\xh(t+1)).
	\end{equation}
	Summing \eqref{eq11} over $t = 1, \ldots, T-1$ and using \eqref{optim}, we obtain
	\begin{align}
	\sum_{t=1}^{T-1} \gt \ip{\xh(t+1)  - \xh(t),\x^\star(t+1)-\xh(t)} &\geq\  \sum_{t=1}^{T-1} \left(\gt-\tfrac{L}{2}\right) \norm{\xh(t+1)-\xh(t)}_2^2  + \tfrac{\mu}{2}\ \sum_{t=1}^{T-1}  \norm{\x^\star(t+1) - \xh(t)}_2^2 \nonumber \\
	&\ \ \ \ {+\sum_{t=1}^{T-1} \ip{\nt,\x^\star(t+1) - \xh(t+1)}}  - \sum_{t=1}^{T-1} \ip{G_t(\xh(t)),\x^\star(t+1) - \xh(t+1)}.
	\end{align}
	{Here, the last two terms can be bounded using the Peter-Paul inequality with parameters $\rho$ and $\eta$ as 
		\begin{align}\label{pprho}
		2\ip{\nt, \xh(t+1)-\x^\star(t+1) } \geq -\tfrac{1}{\rho} \norm{\nt}^2_2 - \rho \norm{\x^\star(t+1) - \xh(t+1)}^2_2,
		\end{align}}
	and
	\begin{align}\label{pp}
	2\ip{G_t(\xh(t)), \xh(t+1)-\x^\star(t+1) } \geq -\tfrac{1}{\eta} \norm{G_t(\xh(t))}^2_2 - \eta\norm{\x^\star(t+1) - \xh(t+1)}^2_2,
	\end{align}
	for $1\leq t \leq T-1$. Using \eqref{pprho}, \eqref{pp} and the definitions of $\gmax$ and $\gmin$, we obtain
	\begin{align}\label{int4}
	\sum_{t=1}^{T-1} \ \ip{\xh(t+1)  - \xh(t),\x^\star(t+1)-\xh(t)} 
	&\geq  \ \tfrac{(2\gmin-L )}{2\gmax}\ \sum_{t=1}^{T-1} \norm{\xh(t+1)-\xh(t)}_2^2   - \tfrac{1}{2\eta\gmax} \Gt - \tfrac{\eta + \rho}{2\gmax} \sum_{t=1}^{T-1} \norm{\x^\star(t+1) - \xh(t+1)}_2^2\nonumber \\
	&\ \ \ { -\tfrac{1}{2\rho\gmax}\sum_{t=1}^{T-1}\norm{\nt}_2^2 }  + \tfrac{\mu}{2\gmax}\ \sum_{t=1}^{T-1}  \norm{\x^\star(t+1) - \xh(t)}_2^2.
	\end{align}
	Recall that $\Gt = \sum_{t=1}^{T-1} \underset{\x\in\X}{\max}\norm{G_t(\x)}_2^2 $. Next we use the expansion
	\begin{align}
	\sum_{t=1}^{T-1}  \norm{\xh(t+1) - \x^\star(t+1)}_2^2  & = \sum_{t=1}^{T-1} \norm{\xh(t+1) - \xh(t)}_2^2 + \sum_{t=1}^{T-1} \norm{\x^\star(t+1)-\xh(t)}_2^2 -2  \sum_{t=1}^{T-1} \ip{\xh(t+1) - \xh(t),\x^\star(t+1)-\xh(t)} \\
	& \leq \left(1 - \tfrac{2\gmin-L}{\gmax}\right) \ \sum_{t=1}^{T-1} \norm{\xh(t+1)-\xh(t)}_2^2
 + (1 -\tfrac{\mu}{\gmax}) \ \sum_{t=1}^{T-1}  \norm{\x^\star(t+1) - \xh(t)}_2^2 { + \tfrac{1}{\rho\gmax}\sum_{t=1}^{T-1}\norm{\nt}_2^2 }  \nonumber \\
	& \ \ \ + \tfrac{1}{\eta\gmax} \Gt + \tfrac{\eta+\rho}{\gmax} \sum_{t=1}^{T-1} \norm{\x^\star(t+1) - \xh(t+1)}_2^2  \label{int3}, 
	\end{align}
	where the last inequality follows from \eqref{int4}. From the Peter-Paul inequality with parameter $\omega$, we have that 
	\begin{align}\label{in3}
	\norm{\x^\star(t+1) - \xh(t) }_2^2\ \leq& (1+\omega ) \norm{\xh(t)-\x^\star(t)  }_2^2  +  (1+\tfrac{1}{\omega})\ \norm{\x^\star(t+1) - \x^\star(t)}_2^2,
	\end{align} 
	for $1\leq t \leq T-1$. Summing over $t=1,\cdots,T-1$ and using the upper bound shown in \eqref{in3} into \eqref{int3} yields, 
	\begin{align}
	(1-\tfrac{\eta+\rho}{\gmax}) \sum_{t=1}^{T-1} \norm{\xh(t+1) - \x^\star(t+1)}_2^2 & \leq \left(\tfrac{\gmax - 2\gmin+L}{\gmax}\right)  \sum_{t=1}^{T-1} \norm{\xh(t+1)-\xh(t)}_2^2  {+ \tfrac{1}{\rho\gmax}\sum_{t=1}^{T-1}\norm{\nt}_2^2 }   \nonumber\\
	&\ \ \  +(1-\tfrac{\mu}{\gmax}) (1+\tfrac{1}{\omega}) \sum_{t=1}^{T-1} \norm{\x^\star(t+1) - \x^\star(t)}_2^2  \nonumber \\
	&\ \ \  + (1 -\tfrac{\mu}{\gmax})  (1+\omega)   \sum_{t=1}^{T-1} \norm{\x^\star(t)-\xh(t)}_2^2  + \tfrac{1}{\eta\gmax} \Gt, 
	\end{align}
	which holds true when $\gmax > \eta+\rho$. Let us further assume that $\eta$ and $\omega$ are chosen such that  $(1-\frac{\mu}{\gmax})(1+\omega) \in  (0, 1-\tfrac{\eta+\rho}{\gmax})$. Recalling that $\gmax > L$ from the statement of Theorem \ref{thm} and that $\xh(1) = \x^\star(1) = \mathbf{s}$, we obtain the following
	\begin{align}\label{fina}
	\s  \left \|\xh(t)  - \x^\star(t) \right \|_2^2 \leq\  \alpha O_T \ + \ \beta\St \ +\ \delta \Gt {\ +\  \zeta \sum_{t=1}^{T-1}\norm{\nt}_2^2 }- (1 -\tfrac{\mu}{\gmax})  (1+\omega) \norm{\xh(T)  - \x^\star(T)}_2^2, 
	\end{align} 
	where $ \alpha = \tfrac{\gmax - 2\gmin+L}{\mu(1+\omega)-\eta - \rho - \omega\gmax }$, $\beta = \tfrac{(\gmax-\mu)(1+\omega)}{\omega(\mu(1+\omega)-\eta-\rho - \omega \gmax)}$, $\delta = \tfrac{1}{\eta(\mu(1+\omega)-\eta-\rho - \omega \gmax)}$, $\zeta = \tfrac{1}{\rho(\mu(1+\omega)-\eta-\rho - \omega \gmax)} $ and 
	$ O_T  = \sum_{t=1}^{T-1} \norm{\xh(t+1) - \xh(t)}_2^2 $. Note that as long as $\omega < \tfrac{\mu-\eta-\rho}{\gmax-\mu}$ and $ \gmax < 2\gmin -L$, the first and the last terms are non-positive and hence can be dropped. Next, we take expectation and use assumption (A5) to yield
	\begin{align}\label{intuppbound}
	\s  &\ \E[\norm{\xh(t)  - \x^\star(t)}_2^2] \leq \beta\St+ \delta\Gt + {\zeta E_T}.
	\end{align}
	{Recall that $E_T = \sum_{t=1}^{T-1} \varepsilon_t^2$. Finally, the regret bound can be obtained by making use of first order convexity condition of the function $-U_t(t)$, which yields:
		\begin{align}\label{foc}
		U_t(\x^\star(t)) &\leq U_t(\xh(t)) + \nabla U_t(\xh(t))^T(\x^\star(t) - \xh(t)),
		\end{align}
		for $1\leq t \leq T$. Summing over $t = 1, \ldots, T$ and rearranging, we obtain
		\begin{align}
		\textbf{Reg}_T \leq \s \ip{\nabla U_t(\xh(t),\x^\star(t)-\xh(t)} \leq G\  \s \norm{\xh(t) - \x^\star(t)}_2, \label{reg2}
		\end{align}
		where we have used the gradient boundedness condition. Now, taking expectation on both sides yields 
		\begin{equation}
		\E \left[\textbf{Reg}_T\right] \leq G\  \s \E \left[\norm{\xh(t) - \x^\star(t)}_2 \right].
		\end{equation}
		Using Cauchy-Schwarz inequality, we have that
		\begin{equation}\label{int7}
		\sum_{t=1}^{T} \E\left[\norm{\xh(t) - \x^\star(t)}_2\right]  \leq \sqrt{T \sum_{t=1}^{T} \E\left[\norm{\xh(t) - \x^\star(t)}_2^2\right]}.
		\end{equation} 
		which, along with \eqref{reg2} and \eqref{fina}, yields the required bound.} 
\end{IEEEproof}
We remark here that it is always possible to choose $\eta$, $\omega$ and $\rho$ so that the aforementioned conditions required for the regret to be sub-linear, are satisfied. A possible choice, obtained by minimizing $\beta$ with respect to $\omega$, $\delta$ with respect to $\eta$ and $\zeta$ with respect to $\rho$ yields
\begin{align}\label{params}
\rho &= \frac{\mu}{2},\ \ \ \ \ \eta = \frac{\mu-\rho}{2} = \frac{\mu}{4},\nonumber  \\  
\omega &= \sqrt{\frac{\gmax-\eta-\rho}{\gmax-\mu}} - 1 \approx \frac{\mu}{8(\gmax - \mu)}, 
\end{align}
which are all constants. Upon substituting \eqref{params} into $\beta, \delta, \zeta$ yields, $\beta = \tfrac{8(\gmax - \mu)(8\gmax - 7\mu)}{\mu^2}$, $\delta = \tfrac{32}{\mu^2}$ and $\zeta = \tfrac{16}{\mu^2}$. Observe that $\beta, \delta$ and $\zeta$ are all constants and are inversely proportional to the square of the strong convexity parameter $\mu$. This limits the applicability of the proposed algorithm to strongly convex loss functions.
To this end, a discussion on how small the parameter $\mu$ can be while guaranteeing sub-linear regret, is added as a remark in the main manuscript.

\subsection{Regret analysis of the alternate update rule }
Consider a scenario where the utility function is known one step a head. This assumption helps in modifying the IOGA update as follows,
\begin{align}
\xh(t+1) = \px\left(\xh(t) + \tfrac{1}{\gamma(t)} \gradt U_{t+1}(\xh(t))\right). \label{ogd}    
\end{align}

\begin{thm}\label{thm2}
	Under the assumptions \textbf{(A1)-(A4)} and for $ \gmax < 2\gmin - L$, $\gmin >  L$ the sequence of $\left\{\x_t\right\}$ generated by IOGA adheres to the regret bound 
	\begin{align}
	\textbf{Reg}_T = \O{\sqrt{T(S_T + E_T)}}.
	\end{align} 
\end{thm}

\begin{IEEEproof}
	Since $-U_{t+1}$ is $\mu$-convex from Assumption \textbf{(A1)}, we have that  for $1\leq t\leq T-1$ 
	\begin{align}\label{sc1}
	U_{t+1}(\x^\star(t+1)) \leq & U_{t+1}(\xh(t)) - \tfrac{\mu}{2} \norm{\x^\star(t+1) - \xh(t)}_2^2  + \ip{\nabla U_{t+1}(\xh(t)),\x^\star(t+1) - \xh(t)}. 
	\end{align}
	Adding and subtracting $\ip{\nabla U_{t+1}(\xh(t)), \xh(t+1)} $ to the right of \eqref{sc1} and rearranging, we obtain  for $1\leq t\leq T-1$ 
	\begin{align}
	U_{t+1}(\x^\star(t+1)) &+ \tfrac{\mu}{2} \norm{\x^\star(t+1) - \xh(t)}_2^2 \nonumber \\ 
	& \leq \  U_{t+1}(\xh(t)) + \ip{\nabla U_{t+1}(\xh(t)),\x^\star(t+1)-\xh(t+1)} + \ip{\nabla U_{t+1}(\xh(t)),\xh(t+1) - \xh(t)}. \label{utdiff1}
	\end{align}
	The IOGA update in \eqref{ogd} can be written as 
	\begin{align}
	\xh(t+1) = &\arg\max_{\x\in\X} \ip{\gradt U_{t+1}(\xh(t)),\x-\xh(t)}  - \tfrac{\gt}{2}\norm{\x-\xh(t)}^2_2, \label{update1}
	\end{align}
	and the optimality condition for \eqref{update1} implies that
	\begin{align}\label{}
	\ip{\gradt U_{t+1}(\xh(t)) - \gt (\xh(t+1)-\xh(t), \x^\star(t+1)-\xh(t+1)} \leq 0, \label{optcon11}
	\end{align}
	since $\x^\star(t+1)\in\X$. Rearranging, \eqref{optcon11} may be re-written as
	\begin{align}
	\ip{\nabla  U_{t+1}(\xh(t)),\ &\x^\star(t+1)-\xh(t+1)} \nonumber\\
	\leq &\ \ \gt \ip{\xh(t+1)-\xh(t),\x^\star(t+1)-\xh(t)} -\gt\norm{\xh(t+1) - \xh(t)}_2^2  -\ip{\nt,\x^\star(t+1) - \xh(t+1)}, \label{int21}
	\end{align} 
	which upon substituting into \eqref{utdiff1} yields for $1\leq t\leq T-1$ 
	\begin{align}\label{scf1}
	U_{t+1}(\x^\star(t+1)) + \frac{\mu}{2} \norm{\x^\star(t+1) - \xh(t)}_2^2 & \leq\ U_{t+1}(\xh(t)) + \ip{\nabla U_{t+1}(\xh(t)),\xh(t+1) - \xh(t)} \nonumber \\ 
	&\ \ \ + \gt\ip{\xh(t+1)- \xh(t), \x^\star(t+1) - \xh(t)} - \gt \norm{\xh(t+1) - \xh(t)}_2^2 \nonumber \\
	& \ \ \ -\ip{\nt,\x^\star(t+1) - \xh(t+1)}.
	\end{align}
	Since $U_{t+1}$ is $L$-smooth from Assumption \textbf{(A2)}, the quadratic lower bound implies that for $1\leq t\leq T-1$ 
	\begin{align}\label{lp1}
	U_{t+1}(\xh(t+1)) \geq  U_{t+1}(\xh(t)) - \tfrac{L}{2}\norm{\xh(t+1)-\xh(t)}_2^2  + \ip{\nabla U_{t+1}(\xh(t)),\xh(t+1)-\xh(t)}.
	\end{align}
	Combining \eqref{scf1} and \eqref{lp1} we obtain for $1\leq t\leq T-1$ 
	\begin{align}\label{eq111}
	U_{t+1}(\x^\star(t+1)) - U_{t+1}(\xh(t+1)) &\leq \left(\tfrac{L}{2} - \gamma(t) \right)\norm{\xh(t+1) - \xh(t)}_2^2 - \tfrac{\mu}{2}\norm{\x^\star(t+1) - \xh(t)}_2^2  \nonumber \\
	&\ \ \ + \gt \ip{\xh(t+1) - \xh(t), \x^\star(t+1) - \xh(t)} -\ip{\nt,\x^\star(t+1) - \xh(t+1)}.
	\end{align}
	Now taking summation over $t$ for $t = 1,\cdots,T-1$ and since $\{\x^\star_t\}_{t=1}^T$ is a solution to off-line version of the problem, it satisfies
	\begin{align} \label{wrg1}
	\s U_{t+1}(\x^\star(t+1)) \geq \s U_{t+1}(\y(t+1)),
	\end{align}
	for any $\{\y(t)\}$ satisfying the coupling constraints. For instance, the sequence $\{\xh(t) \in \X\}_{t=1}^T$ satisfies $g_t(\x(t),\x(t+1)) \leq 0$ for all $1\leq t \leq T-1$ from assumption (\textbf{A4}). Now by employing \eqref{wrg1} we get 
	\begin{align}
	\s \gt \ \ip{\xh(t+1) - \xh(t), \x^\star(t+1) - \xh(t)} 
	&\geq\  \s \left(\gt-\tfrac{L}{2}\right) \norm{\xh(t+1)-\xh(t)}_2^2  + \tfrac{\mu}{2}\ \s  \norm{\x^\star(t+1) - \xh(t)}_2^2 \nonumber \\
	& \ \ \ +\s  \ip{\nt,\x^\star(t+1) - \xh(t+1)}.
	\end{align}	
	Further we can write 
	\begin{align}\label{int41}
	\s \ \ip{\xh(t+1) - \xh(t),\x^\star(t+1)-\xh(t)} & \geq  \ \tfrac{(2\gmin-L )}{2\gmax}\ \s \norm{\xh(t+1)-\xh(t)}_2^2  + \tfrac{\mu}{2\gmax}\ \s  \norm{\x^\star(t+1) - \xh(t)}_2^2 \nonumber \\
	& \ \ \ +\tfrac{1}{\gmax} \s  \ip{\nt,\x^\star(t+1) - \xh(t+1)},
	\end{align}
	where $\gmax = \underset{t}{\text{arg}\max } \left\{\gt \right\}$ and $\gmin = \underset{t}{\text{arg}\min } \left\{\gt \right\}$. Further, we have the following 
	\begin{align}
	\s  \norm{\xh(t+1) - \x^\star(t+1)}_2^2 & = \s \norm{\xh(t+1) - \xh(t)}_2^2 + \s \norm{\x^\star(t+1)-\xh(t)}_2^2 -2\ \s \ip{\xh(t+1) - \xh(t),x^\star(t+1)-\xh(t)} \\
	&\leq \left(\tfrac{\gmax-2\gmin+L}{\gmax}\right) \ \s \norm{\xh(t+1)-\xh(t)}_2^2
+ (1 -\tfrac{\mu}{\gmax}) \ \s  \norm{\x^\star(t+1) - \xh(t)}_2^2 \nonumber \\
	& \hspace{0.3cm} -\tfrac{2}{\gmax}  \label{int31} \s  \ip{\nt,\x^\star(t+1) - \xh(t+1)},
	\end{align}
	where the last inequality follows from \eqref{int41}. Next from Peter-Paul inequality with parameter $\eta > 0$ we have
	\begin{align}\label{in31}
	\norm{\x^\star(t+1) - \xh(t) }_2^2\ \leq\ (1+\eta )\ \norm{\xh(t)-\x^\star(t)  }_2^2\ +\  (1+\tfrac{1}{\eta})\ \norm{\x^\star(t+1) - \x^\star(t)}_2^2. 
	\end{align} 
	Next using the upper bound shown in \eqref{in31} into \eqref{int31}, we get  
	\begin{align}
    \s \norm{\xh(t+1) - \x^\star(t+1)}_2^2 &\leq \left(\tfrac{\gmax - 2\gmin+L}{\gmax}\right) \s \norm{\xh(t+1)-\xh(t)}_2^2+(1-\tfrac{\mu}{\gmax}) (1+\tfrac{1}{\eta})\s \norm{\x^\star(t+1) - \x^\star(t)}_2^2  \nonumber \\
	&\ \ \  + (1 -\tfrac{\mu}{\gmax})  (1+\eta)  \s \norm{\xh(t)-\x^\star(t)}_2^2 + \tfrac{1}{\rho \gmax} \norm{\nt}^2_2 + \tfrac{\rho}{\gmax} \norm{\x^\star(t+1) - \xh(t+1)}_2^2.
	\end{align}
	Next for $(1-\frac{\mu}{\gmax})(1+\eta) \in  (0, 1)$ while $\gmax > L$,   
	and using the fact that $\xh(1) = \x^\star(1)$, we obtain the following
	\begin{align}\label{fina1}
	\sum_{t=1}^{T} \left \|\xh(t)  - \x^\star(t) \right \|_2^2 \leq \alpha\ O_T + \beta\ \St + \delta\ E_T  - (1 -\tfrac{\mu}{\gmax})  (1+\eta) \norm{\xh(T)  - \x^\star(T)}_2^2,	  
	\end{align} 
	where $ \alpha = \tfrac{\gmax - 2\gmin+L}{\mu(1+\eta)- \eta \gmax - \rho}$, $\beta = \tfrac{(\gmax-\mu)(1+\eta)}{\eta(\mu(1+\eta)-\eta \gmax- \rho) }$, $\delta = \tfrac{1}{\rho(\mu(1+\eta)-\eta \gmax - \rho) }$ 
	$ O_T  = \s \norm{\xh(t+1) - \xh(t)}_2^2 $ and $\St=\sum_{t=1}^{T-1} \norm{\x^\star(t+1) - \x^\star(t)}_2^2 $. 
	Now, observe that first term can be dropped for $ \gmax < 2\gmin -L$, also the last term. Finally, we obtain the required bound by making use of first order convexity condition of the function $-U_t(t)$, Cauchy-Schwartz inequality and boundedness of gradient as follows for $1\leq t \leq T$ 
	\begin{equation}\label{foc1}
	U_t(\x^\star(t))\leq U_t(\xh(t)) + \nabla U_t(\xh(t))^T(\x^\star(t) - \xh(t)).
	\end{equation} 	 
	Taking summation over $t$ on both sides we get 
	\begin{align}\label{int61}
	\sum_{t=1}^{T} U_t(\x^\star(t)) -U_t(\xh(t))  & \leq\ \s \nabla U_t(\xh(t))^\intercal (\x^\star(t)-\xh(t)) \nonumber \\
	& \leq G\  \s \norm{\xh(t) - \x^\star(t)}_2. 
	\end{align} 	 
	It holds in general that 
	\begin{equation}\label{int71}
	\sum_{t=1}^{T} \norm{\xh(t) - \x^\star(t)}_2  \leq \sqrt{T \sum_{t=1}^{T} \norm{\xh(t) - \x^\star(t)}_2^2 }.
	\end{equation} 	  
	Now utilizing \eqref{fina1}, \eqref{int71} into \eqref{int61} we get the following
	\begin{equation}
	\sum_{t=1}^{T} U_t(\x^\star(t)) - U_t(\xh(t))  \leq \O{\sqrt{T (\St + E_T)}}.
	\end{equation} 	
\end{IEEEproof}

 \end{document}